\documentclass{article}
\usepackage{latexsym}
\usepackage[varumlaut]{yfonts}
\usepackage{epsfig}
\usepackage{amsmath}
\usepackage{amssymb}
\usepackage{amsthm}
\usepackage{eufrak}
\newtheorem{theorem}{Theorem}[section]
\newtheorem{corollary}[theorem]{Corollary}

\newtheorem{lemma}[theorem]{Lemma}
\newtheorem{conjecture}[theorem]{Conjecture}

\theoremstyle{remark}
\newtheorem{remark}[theorem]{Remark}

\theoremstyle{definition}
\newtheorem{definition}[theorem]{Definition}

\DeclareMathOperator\diag{diag}

\DeclareMathOperator\rk{rk}

\DeclareMathOperator\sign{sgn}

\DeclareMathOperator\tr{tr}
\DeclareMathOperator\spa{span}
\DeclareMathOperator\Supp{supp}

\begin{document}

\title{Copositive matrices with circulant zero support set}

\author{Roland Hildebrand \thanks{%
LJK, 700 avenue Centrale, Domaine Universitaire, 38401 St. Martin d'H\`eres, France ({\tt roland.hildebrand@imag.fr}).}}

\maketitle

\begin{abstract}
Let $n \geq 5$ and let $u^1,\dots,u^n$ be nonnegative real $n$-vectors such that the indices of their positive elements form the sets $\{1,2,\dots,n-2\},\{2,3,\dots,n-1\},\dots,\{n,1,\dots,n-3\}$, respectively. Here each index set is obtained from the previous one by a circular shift. The set of copositive forms which vanish on the vectors $u^1,\dots,u^n$ is a face of the copositive cone ${\cal C}^n$. We give an explicit semi-definite description of this face and of its subface consisting of positive semi-definite forms, and study their properties. If the vectors $u^1,\dots,u^n$ and their positive multiples exhaust the zero set of an exceptional copositive form belonging to this face, then we say it has minimal circulant zero support set, and otherwise non-minimal circulant zero support set. We show that forms with non-minimal circulant zero support set are always extremal, and forms with minimal circulant zero support sets can be extremal only if $n$ is odd. We construct explicit examples of extremal forms with non-minimal circulant zero support set for any order $n \geq 5$, and examples of extremal forms with minimal circulant zero support set for any odd order $n \geq 5$. The set of all forms with non-minimal circulant zero support set, i.e., defined by different collections $u^1,\dots,u^n$ of zeros, is a submanifold of codimension $2n$, the set of all forms with minimal circulant zero support set a submanifold of codimension $n$.
\end{abstract}




\section{Introduction}

Let ${\cal S}^n$ be the vector space of real symmetric $n \times n$ matrices. In this space, we may define the cone ${\cal N}^n$ of element-wise nonnegative matrices, the cone ${\cal S}_+^n$ of positive semi-definite matrices, and the cone ${\cal C}^n$ of copositive matrices, i.e., matrices $A \in {\cal S}^n$ such that $x^TAx \geq 0$ for all $x \in \mathbb R_+^n$. Obviously we have ${\cal S}_+^n + {\cal N}^n \subset {\cal C}^n$, but the converse inclusion holds only for $n \leq 4$ \cite[Theorem 2]{Diananda62}. Copositive matrices which are not elements of the sum ${\cal S}_+^n + {\cal N}^n$ are called {\it exceptional}. Copositive matrices play an important role in non-convex and combinatorial optimization, see, e.g., \cite{Burer09} or the surveys \cite{Duer10},\cite{HUS10},\cite{BSU12},\cite{Dickinson_2013}. Of particular interest are the exceptional extreme rays of ${\cal C}^n$.

A fruitful concept in the study of copositive matrices is that of zeros and their supports, initiated in the works of Baumert \cite{Baumert66},\cite{Baumert67}, see also \cite{Hildebrand14a},\cite{DickinsonHildebrand16}, and \cite{ShakedMonderer15} for further developments and applications. A non-zero vector $u \in \mathbb R_+^n$ is called a {\it zero} of a copositive matrix $A \in {\cal C}^n$ if $u^TAu = 0$. The {\it support} $\Supp u$ of a zero $u$ is the index set of its positive elements.

Note that each of the cones ${\cal N}^n$, ${\cal S}_+^n$, ${\cal C}^n$ is invariant with respect to a simultaneous permutation of the row and column indices, and with respect to a simultaneous pre- and post-multiplication with a positive definite diagonal matrix. These operations generate a group of linear transformations of ${\cal S}^n$, which we shall call ${\cal G}_n$.

Exceptional copositive matrices first appear at order $n = 5$. The {\it Horn matrix}
\begin{equation} \label{Horn_matrix}   H = \left( \begin{array}{rrrrr}
              1 & -1 & 1 & 1 & -1 \\
              -1 & 1 & -1 & 1 & 1 \\
              1 & -1 & 1 & -1 & 1 \\
              1 & 1 & -1 & 1 & -1 \\
              -1 & 1 & 1 & -1 & 1
      \end{array} \right),
\end{equation}
named after its discoverer Alfred Horn, and the other matrices in its ${\cal G}_5$-orbit have been the first examples of exceptional extremal copositive matrices \cite{HallNewman63}. Any other exceptional extremal matrix in ${\cal C}^5$ lies in the ${\cal G}_5$-orbit of a matrix
\begin{equation} \label{T_matrices}
T(\theta) = \begin{pmatrix}
       1 & -\cos\theta_1 & \cos(\theta_1+\theta_2) & \cos(\theta_4+\theta_5) & -\cos\theta_5 \\
       -\cos\theta_1 & 1 & -\cos\theta_2 & \cos(\theta_2+\theta_3) & \cos(\theta_5+\theta_1) \\
       \cos(\theta_1+\theta_2) & -\cos\theta_2 & 1 & -\cos\theta_3 & \cos(\theta_3+\theta_4) \\
       \cos(\theta_4+\theta_5) & \cos(\theta_2+\theta_3) & -\cos\theta_3 & 1 & -\cos\theta_4 \\
        -\cos\theta_5 & \cos(\theta_5+\theta_1) & \cos(\theta_3+\theta_4) & -\cos\theta_4 & 1
   \end{pmatrix}
\end{equation}
for some angles $\theta_k \in (0,\pi)$ satisfying $\sum_{k=1}^5 \theta_k < \pi$ \cite[Theorem 3.1]{Hildebrand12a}. Both the Horn matrix $H$ and the matrices $T(\theta)$, which are referred to as Hildebrand matrices, possess zeros with supports $\{1,2,3\},\{2,3,4\},\{3,4,5\},\{4,5,1\},\{5,1,2\}$, respectively. Note that these supports are exactly the vertex subsets obtained by removing the vertices of a single edge in the cycle graph $C_5$.

In this contribution we shall generalize the exceptional extremal elements of ${\cal C}^5$ to arbitrary order $n \geq 5$ by taking the above property of the supports as our point of departure. Fix a set ${\bf u} = \{ u^1,\dots,u^n \} \subset \mathbb R_+^n$ of nonnegative vectors with supports $\{1,2,\dots,n-2\}$, $\{2,3,\dots,n-1\},\dots,\{n,1,\dots,n-3\}$, respectively, i.e., the supports of the vectors $u^j$ are the vertex subsets obtained by removing the vertices of a single edge in the cycle graph $C_n$. We then consider the faces
\begin{equation} \label{faces_def}
F_{\bf u} = \{ A \in {\cal C}^n \,|\, (u^j)^TAu^j = 0\ \forall\ j = 1,\dots,n \}, \qquad P_{\bf u} = \{ A \in {\cal S}_+^n \,|\, (u^j)^TAu^j = 0\ \forall\ j = 1,\dots,n \}
\end{equation}
of the copositive cone and the positive semi-definite cone, respectively. Note that $P_{\bf u} \subset F_{\bf u}$.

One of our main results is an explicit semi-definite description of the faces $F_{\bf u}$ and $P_{\bf u}$ (Theorem \ref{thm:sdp_rep}). In order to obtain this description, we associate the set ${\bf u}$ with a discrete-time linear dynamical system ${\bf S}_{\bf u}$ of order $d = n-3$ and with time-dependent coefficients having period $n$. If ${\cal L}_{\bf u}$ is the $d$-dimensional solution space of this system, then there exists a canonical bijective linear map between $F_{\bf u}$ and the set of positive semi-definite symmetric bilinear forms on the dual space ${\cal L}_{\bf u}^*$ satisfying certain additional homogeneous linear equalities and inequalities. For an arbitrary collection ${\bf u}$ in general only the zero form satisfies the corresponding linear matrix inequality (LMI) and the face $F_{\bf u}$ consists of the zero matrix only. However, for every $n \geq 5$ there exist collections ${\bf u}$ for which the LMI has non-trivial feasible sets.

The properties of the copositive matrices in $F_{\bf u}$ are closely linked to the properties of the periodic linear dynamical system ${\bf S}_{\bf u}$. Such systems are the subject of {\it Floquet theory}, see, e.g., \cite[Section 3.4]{Elaydi}. We need only the concept of the {\it monodromy operator} and its eigenvalues, the {\it Floquet multipliers}, which we shall review in Section \ref{Floquet}. We show that the face $P_{\bf u}$ is isomorphic to ${\cal S}_+^{d_1}$, where $d_1$ is the geometric multiplicity of the Floquet multiplier 1, or equivalently, the dimension of the subspace of $n$-periodic solutions of ${\bf S}_{\bf u}$. For the existence of exceptional copositive matrices in $F_{\bf u}$ it is necessary that all or all but one Floquet multiplier are located on the unit circle (Corollary \ref{cor:corank1}).

We are able to describe the structure of $F_{\bf u}$ explicitly in general. Exceptional matrices $A \in F_{\bf u}$ can be divided in two categories. If every zero of $A$ is proportional to one of the zeros $u^1,\dots,u^n$, then we say the copositive matrix $A$ has {\it minimal circulant zero support set}, otherwise it has {\it non-minimal circulant zero support set}. We show that matrices with non-minimal circulant zero support set are always extremal, while matrices with minimal circulant zero support set can be extremal only for odd $n$. For even $n$ a matrix with minimal circulant zero support set can be represented as a non-trivial sum of a matrix with non-minimal circulant zero support set and a positive semi-definite rank 1 matrix, the corresponding face $F_{\bf u}$ is then isomorphic to $\mathbb R_+^2$. For odd $n$ a sufficient condition for extremality of a matrix with minimal circulant zero support set is that $-1$ does not appear among the Floquet multipliers (Theorem \ref{thm:reg_deg}). For every $n \geq 5$ the matrices with non-minimal circulant zero support set constitute an algebraic submanifold of ${\cal S}^n$ of codimension $2n$ (Theorem \ref{thm:manifold_2n}), while the matrices with minimal circulant zero support set form an algebraic submanifold of codimension $n$ (Theorem \ref{thm:manifold_n}), in which the extremal matrices form an open subset (Theorem \ref{thm:manifold_n_ext}).

Finally, in Section \ref{subs:circulant} we construct explicit examples of circulant (i.e., invariant with respect to simultaneous circular shifts of row and column indices) exceptional extremal copositive matrices, both with minimal and non-minimal circulant zero support set. We also give an exhaustive description of all exceptional matrices for order $n = 6$ with non-minimal circulant zero support set in Section \ref{sec:n6}. Some auxiliary results whose proofs would interrupt the flow of exposition are collected in two appendices.

\subsection{Further notations}

For $n \geq 5$ an integer, define the ordered index sets $I_1 = (1,2,\dots,n-2)$, $I_2 = (2,3,\dots,n-1),\dots,I_n = (n,1,\dots,n-3)$ of cardinality $n-2$, each obtained by a circular shift of the indices from the previous one. We will need also the index sets $I_1' = (1,2,\dots,n-3),\dots,I_n' = (n,1,\dots,n-4)$ defined similarly.

Let $k > 0$ be an integer. For a vector $u \in \mathbb R^k$, a $k \times k$ matrix $M$, and an ordered index set $I \subset \{1,\dots,k\}$ of cardinality $|I|$, we shall denote by $u_i$ the $i$-th entry of $u$, by $u_I$ the subvector $(u_i)_{i \in I} \in \mathbb R^{|I|}$ of $u$ composed of the elements with index in the ordered set $I$, by $M_{ij}$ the $(i,j)$-th entry of $M$, and by $M_I$ the principal submatrix $(M_{ij})_{i,j \in I} \in \mathbb R^{|I| \times |I|}$ of $M$ composed of the elements having row and column index in $I$.

In order to distinguish it from the index sets $I_1,\dots,I_n$ defined above, we shall denote the identity matrix or the identity operator by $\it{Id}$ or $\it{Id}_k$ if it is necessary to indicate the size of the matrix. Denote by $E_{ij} \in {\cal N}^n$ the matrix which has ones at the positions $(i,j)$ and $(j,i)$ and zeros elsewhere. For a real number $r$, we denote by $\lfloor r \rfloor$ the largest integer not exceeding $r$ and by $\lceil r \rceil$ the smallest integer not smaller than $r$.

{\definition \label{def:ref} Let $A \in {\cal C}^n$ be an exceptional copositive matrix possessing zeros $u^1,\dots,u^n \in \mathbb R_+^n$ such that $\Supp u^j = I_j$, $j = 1,\dots,n$. We say the matrix $A$ has {\sl minimal circulant zero support set} if every zero $v$ of $A$ is proportional to one of the zeros $u^1,\dots,u^n$, and {\sl non-minimal circulant zero support set} otherwise. }

\medskip

A real symmetric $n \times n$ matrix will be called \emph{positive semi-definite} (PSD), denoted $A \succeq 0$, if $x^TAx \geq 0$ for all $x \in \mathbb R^n$. A vector $u \in \mathbb R^n$ will be called \emph{nonnegative}, denoted $u \geq 0$, if $u_i \geq 0$ for all $i = 1,\dots,n$.

\section{Conditions for copositivity}

In this section we consider matrices $A \in {\cal S}^n$ such that the submatrices $A_{I_1},\dots,A_{I_n}$ are all positive semi-definite and possess element-wise positive kernel vectors. We derive necessary and sufficient conditions for such a matrix to be copositive. The goal of the section is to prove Theorem \ref{first_cond} below. We start with a few simple auxiliary lemmas.

%


{\lemma \cite[Lemma 2.4]{DDGH13a} \label{zeroPSD} Let $A \in {\cal C}^n$ and let $u$ be a zero of $A$. Then the principal submatrix $A_{\Supp{u}}$ is positive semi-definite. }


{\lemma \cite[p.200]{Baumert66} \label{first_order} Let $A \in {\cal C}^n$ and let $u$ be a zero of $A$. Then $Au \geq 0$. }

%
%
%
%
%
%
%

{\lemma \label{lem:irred_N} Let $n \geq 5$, and let $i,j \in \{1,\dots,n\}$ be arbitrary indices. Then there exists $k \in \{1,\dots,n\}$ such that $i,j \in I_k$. }

\begin{proof}
For every index $i \in \{1,\dots,n\}$, there exist exactly two indices $k$ such that $i \not\in I_k$. The assertion of the lemma then follows from the Dirichlet principle.
\end{proof}

{\corollary \label{cor:irred_N} Let $n \geq 5$ and let $A \in {\cal C}^n$ have zeros $u^1,\dots,u^n$ with supports $I_1,\dots,I_n$, respectively. Then for every pair of indices $i,j \in \{1,\dots,n\}$, the matrix $A - \varepsilon E_{ij}$ is not copositive for every $\varepsilon > 0$. }

\begin{proof}
Let $i,j \in \{1,\dots,n\}$. By Lemma \ref{lem:irred_N} there exists a zero $u^k$ of $A$ such that $i,j \in \Supp u^k$. The assertion then follows by $(u^k)^T(A - \varepsilon E_{ij})u^k < 0$ for every $\varepsilon > 0$, see \cite[p. 10]{Baumert67}.
\end{proof}

{\corollary \cite[Corollary 3.7]{Baumert66} \label{Baumert66Cor37} Let $A \in {\cal C}^n$ be such that $A - \varepsilon E_{ij}$ is not copositive for every $\varepsilon > 0$ and all $i,j = 1,\dots,n$. If $A$ has a zero with support of cardinality $n-1$, then $A$ is positive semi-definite. }

\medskip

A zero $u \in \mathbb R_+^n$ of $A$ is called {\it minimal} if there is no other zero $v \in \mathbb R_+^n$ of $A$ such that the inclusion $\Supp v \subset \Supp u$ is strict. Minimality is hence a property of the support rather than of the zero itself.

{\lemma \label{minimal_corank1} \cite[Lemma 3.7]{Hildebrand14a} Let $A \in {\cal C}^n$ and let $I \subset \{1,\dots,n\}$ be an index set. Then $I$ is the support of a minimal zero of $A$ if and only if the principal submatrix $A_I$ is positive semi-definite of corank 1 and has a kernel vector with positive elements. }

For a fixed copositive matrix, the number of its minimal zeros whose elements sum up to 1 is finite \cite[Corollary 3.6]{Hildebrand14a}. The next result furnishes a criterion involving minimal zeros to check whether a copositive matrix is extremal.

{\lemma \label{extreme_criterion} \cite[Theorem 17]{DickinsonHildebrand16} Let $A \in {\cal C}^n$, and let $u^1,\dots,u^m$ be its minimal zeros whose elements sum up to 1. Then $A$ is extremal if and only if the solution space of the linear system of equations in the matrix $X \in {\cal S}^n$ given by
\[ (Xu^j)_i = 0\qquad\forall\ i,j:\ (Au^j)_i = 0
\]
is one-dimensional (and hence generated by $A$). }

{\lemma \label{irred_rank1} \cite[Lemma 4.3]{Hildebrand14a} Let $A \in {\cal C}^n$ be a copositive matrix and let $w \in \mathbb R^n$. Then there exists $\varepsilon > 0$ such that $A - \varepsilon ww^T$ is copositive if and only if $\langle w,u \rangle = 0$ for all zeros $u$ of $A$. }

\bigskip

These results allow us to prove the following theorem on matrices $A \in {\cal S}^n$ having zeros $u^1,\dots,u^n$ with supports $I_1,\dots,I_n$, respectively.

{\theorem \label{first_cond} Let $n \geq 5$ and let $A \in {\cal S}^n$ be such that for every $j = 1,\dots,n$ there exists a nonnegative vector $u^j$ with $\Supp u^j = I_j$ satisfying $(u^j)^TAu^j = 0$. Then the following are equivalent:
\begin{itemize}
\setlength{\itemsep}{1pt}
\setlength{\parskip}{0pt}
\setlength{\parsep}{0pt}
\item[(i)] $A$ is copositive;
\item[(ii)] every principal submatrix of $A$ of size $n-1$ is copositive;
\item[(iii)] every principal submatrix of $A$ of size $n-1$ is in ${\cal S}_+^{n-1} + {\cal N}^{n-1}$;
\item[(iv)] $A_{I_j}$ is positive semi-definite for $j = 1,\dots,n$, $(u^n)^TAu^1 \geq 0$, and $(u^j)^TAu^{j+1} \geq 0$ for $j = 1,\dots,n-1$.
\end{itemize}

Moreover, given above conditions (i)---(iv), the following are equivalent:
\begin{itemize}
\setlength{\itemsep}{1pt}
\setlength{\parskip}{0pt}
\setlength{\parsep}{0pt}
\item[(a)] $A$ is positive semi-definite;
\item[(b)] at least one of the $n$ numbers $(u^n)^TAu^1$ and $(u^j)^TAu^{j+1}$, $j = 1,\dots,n-1$, is zero;
\item[(c)] all $n$ numbers $(u^n)^TAu^1$ and $(u^j)^TAu^{j+1}$, $j = 1,\dots,n-1$, are zero;
\item[(d)] $A$ is not exceptional.
\end{itemize}}

\begin{proof}
(i) $\Rightarrow$ (iv) is a consequence of Lemmas \ref{zeroPSD} and \ref{first_order}.

(iv) $\Rightarrow$ (iii) is a consequence of Lemma \ref{completion} in the Appendix, applied to the $(n-1) \times (n-1)$ principal submatrices of $A$.

(iii) $\Rightarrow$ (ii) is trivial.

(ii) $\Rightarrow$ (i): Let $\Delta = \{ v = (v_1,\dots,v_n)^T \in \mathbb R_+^n \,|\, \sum_{j=1}^n v_j = 1 \}$ be the standard simplex. By (ii) the quadratic form $A$ is nonnegative on $\partial\Delta$. On $\partial\Delta$ it then reaches its global minimum 0 at appropriate positive multiples $\alpha_ju^j \in \Delta$ of the zeros $u^j$ for all $j = 1,\dots,n$. Since the line segment connecting $\alpha_1u^1$ and $\alpha_2u^2$ still lies in $\partial\Delta$, the quadratic function $Q(v) = v^TAv$ cannot be strictly convex on $\Delta$. But then it reaches its global minimum over $\Delta$ on the boundary $\partial\Delta$. This minimum over $\Delta$ then also equals 0, which proves the copositivity of $A$.

We have shown the equivalence of conditions (i)---(iv). Let us now assume that $A$ satisfies (i)---(iv) and pass to the second part of the theorem.

(a) $\Rightarrow$ (c): If $A \succeq 0$, then all vectors $u^j$, $j = 1,\dots,n$, are in the kernel of $A$. This implies (c).

(c) $\Rightarrow$ (b) is trivial.

(b) $\Rightarrow$ (a): Without loss of generality, let $(u^1)^TAu^2 = 0$. It then follows that $u_+ = u^1 + u^2$ is also a zero of $A$, with $\Supp u_+ = \{1,\dots,n-1\}$. By Corollaries \ref{cor:irred_N} and \ref{Baumert66Cor37} $A$ is then positive semi-definite, which proves (a).

(a) $\Rightarrow$ (d) holds by definition.

(d) $\Rightarrow$ (a): Assume that $A$ can be written as a sum $A = P + N$ with $P \in {\cal S}_+^n$ and $N \in {\cal N}^n$. By Corollary \ref{cor:irred_N} none of the elements of $N$ can be positive and hence $A = P$, implying (a).
\end{proof}

Let us comment on Theorem \ref{first_cond}. It states that the presence of $n$ zeros with supports $I_j$, $j = 1,\dots,n$ places stringent constraints on a copositive matrix $A \in {\cal C}^n$. Such a matrix must either be exceptional or positive semi-definite. Which of these two cases arises is determined by any of the $n$ numbers in condition (iv) of the theorem, which are either simultaneously positive or simultaneously zero.

\section{Linear systems with periodic coefficients} \label{Floquet}

In this section we investigate the solution spaces of linear periodic dynamical systems and perform some linear algebraic constructions on them. These will be later put in correspondence to copositive forms. First we shall introduce the monodromy and the Floquet multipliers associated with such systems, for further reading about these and related concepts see, e.g., \cite[Section 3.4]{Elaydi}.

We consider real scalar discrete-time homogeneous linear dynamical systems governed by the equation
\begin{equation} \label{periodic_system}
x_{t+d} + \sum_{i=0}^{d-1} c^t_ix_{t+i} = \sum_{i=0}^d c^t_ix_{t+i} = 0,\qquad t = 1,2,\dots
\end{equation}
where $x_t \in \mathbb R$ is the value of the solution $x$ at time instant $t$, $d > 0$ is the order, and $c^t = (c^t_0,\dots,c^t_d)^T \in \mathbb R^{d+1}$, $t \geq 1$, are the coefficient vectors of the system. For convenience we have set $c^t_d = 1$ for all $t \geq 1$. We assume that the coefficients are periodic with period $n > d$, i.e., $c^{t+n} = c^t$ for all $t \geq 1$. Denote by ${\cal L}$ the linear space of all solutions $x = (x_t)_{t \geq 1}$. This space has dimension $d$ and can be parameterized, e.g., by the vector $(x_1,\dots,x_d) \in \mathbb R^d$ of initial conditions.

If $x = (x_t)_{t \geq 1}$ is a solution of the system, then $y = (x_{t+n})_{t \geq 1}$ is also a solution by the periodicity of the coefficients. The corresponding linear map ${\mathfrak M}: {\cal L} \to {\cal L}$ taking $x$ to $y$ is called the {\it monodromy} of the periodic system. Its eigenvalues are called {\it Floquet multipliers}. The following result is a trivial consequence of this definition.

{\lemma \label{periodic1} Let ${\cal L}_{per} \subset {\cal L}$ be the subspace of $n$-periodic solutions of system \eqref{periodic_system}. Then $x \in {\cal L}_{per}$ if and only if $x$ is an eigenvector of the monodromy operator $\mathfrak M$ with eigenvalue 1. In particular, $\dim {\cal L}_{per}$ equals the geometric multiplicity of the eigenvalue 1 of $\mathfrak M$. \qed }

We call system \eqref{periodic_system} \emph{time-reversible} if for every $T > 0$, and for every solution $(x_t)_{t \geq 1}$, the values $(x_t)_{t < T}$ can be restored from the values $(x_t)_{t \geq T}$. Clearly this is the case if and only if $c^t_0 \not= 0$ for all $t = 1,\dots,n$, in this case we may express $x_t$ as a function of $x_{t+1},\dots,x_{t+d}$. In fact, the following result holds.

{\lemma \label{monodromy_det} The determinant of the monodromy operator is given by $\det {\mathfrak M} = (-1)^{nd} \prod_{t=1}^n c^t_0$. }

\begin{proof}
From \eqref{periodic_system} it follows that the determinant of the linear map taking the vector $(x_t,x_{t+1},\dots,x_{t+d-1})$ to $(x_{t+1},\dots,x_{t+d})$ equals $(-1)^dc^t_0$. Iterating this map for $t = 1,\dots,n$, we get that the determinant of the linear map taking the vector $(x_1,\dots,x_d)$ to $(x_{n+1},\dots,x_{n+d})$ equals $(-1)^{nd} \prod_{t=1}^n c^t_0$. The claim now follows from the fact that the vector $(x_1,\dots,x_d)$ parameterizes the solution space ${\cal L}$.
\end{proof}

Let us now consider the space ${\cal L}^*$ of linear functionals on the solution space ${\cal L}$. For every $t \geq 1$, the map taking a solution $x = (x_s)_{s \geq 1}$ to its value $x_t$ at time instant $t$ is such a linear functional. We shall denote this evaluation functional by ${\bf e}_t \in {\cal L}^*$. By definition of the monodromy we have ${\bf e}_{t+n} = {\mathfrak M}^*{\bf e}_t$ for all $t \geq 1$, where ${\mathfrak M}^*: {\cal L}^* \to {\cal L}^*$ is the adjoint of $\mathfrak M$. Moreover,
\begin{equation} \label{e_relation}
\sum_{i=0}^d c^t_i{\bf e}_{t+i} = 0\qquad \forall\ t \geq 1
\end{equation}
as a consequence of \eqref{periodic_system}.

Our main tool in the study of copositive forms in this paper are positive semi-definite symmetric bilinear forms $B$ on ${\cal L}^*$ which are invariant with respect to a time shift by the period $n$, i.e.,
\begin{equation} \label{periodic_Q}
B({\bf e}_{t+n},{\bf e}_{s+n}) = B({\bf e}_t,{\bf e}_s)\qquad \forall\ t,s \geq 1.
\end{equation}

{\definition We call a symmetric bilinear form $B$ on ${\cal L}^*$ satisfying relation \eqref{periodic_Q} a {\it shift-invariant} form. }

{\lemma \label{lem:inv} Assume above notations. A symmetric bilinear form $B$ on ${\cal L}^*$ is shift-invariant if and only if it is preserved by the adjoint of the monodromy, i.e., $B(w,w') = B({\mathfrak M}^*w,{\mathfrak M}^*w')$ for all $w,w' \in {\cal L}^*$. An equivalent set of conditions is given by $B({\bf e}_{t+n},{\bf e}_{s+n}) = B({\bf e}_t,{\bf e}_s)$ for all $t,s \in \{1,\dots,d\}$. }

\begin{proof}
The assertions hold because the evaluation functionals ${\bf e}_1,\dots,{\bf e}_d$ form a basis of ${\cal L}^*$ and ${\bf e}_{t+n} = {\mathfrak M}^*{\bf e}_t$ for all $t \geq 1$.
\end{proof}

The shift-invariant forms are hence determined by a finite number of linear homogeneous equations and constitute a linear subspace of the space of symmetric bilinear forms on ${\cal L}^*$.

The space of symmetric bilinear forms on ${\cal L}^*$ can be viewed as the space of symmetric contra-variant second order tensors over ${\cal L}$, i.e., it is the linear hull of tensor products of the form $x \otimes x$, $x \in {\cal L}$. It is well-known that a symmetric bilinear form can be diagonalized, i.e., represented as a finite sum $B = \sum_{k=1}^r \sigma_k x^k \otimes x^k$ with $\sigma_k \in \{-1,+1\}$, $x^k \in {\cal L}$, $k = 1,\dots,r$, the vectors $x^k$ being linearly independent. The vectors $x^k$ in this decomposition are not unique, but their number $r$ and their linear hull depend only on $B$ and are called the {\it rank} $\rk\,B$ and the {\it image} $\it{Im}\,B$ of $B$, respectively. The form is positive semi-definite if all coefficients $\sigma_k$ in its decomposition equal 1.

{\lemma \label{lem:inv_unitary} Let $B$ be a shift-invariant symmetric positive semi-definite bilinear form $B$ on ${\cal L}^*$, of rank $r$. Then there exist at least $r$ (possibly complex) linearly independent eigenvectors of $\mathfrak M$ with eigenvalues on the unit circle. }

\begin{proof}
Let $B = \sum_{k=1}^r x^k \otimes x^k$ be a decomposition of $B$ as above and complete the linearly independent set $\{x^1,\dots,x^r\}$ to a basis $\{x^1,\dots,x^d\}$ of ${\cal L}$. In the coordinates defined by this basis and its dual basis in ${\cal L}^*$ the form $B$ is then given by the diagonal matrix $\diag(\it{Id}_r,0,\dots,0)$. Let $M$ be the coefficient matrix of $\mathfrak M$ in this basis, partitioned into submatrices $M_{11},M_{12},M_{21},M_{22}$ corresponding to the partition of the basis into subsets $\{x^1,\dots,x^r\}$ and $\{x^{r+1},\dots,x^d\}$. Then by Lemma \ref{lem:inv} the shift-invariance of $B$ is equivalent to the condition
\[ \begin{pmatrix} \it{Id}_r & 0 \\ 0 & 0 \end{pmatrix} = \begin{pmatrix} M_{11} & M_{12} \\ M_{21} & M_{22} \end{pmatrix} \begin{pmatrix} \it{Id}_r & 0 \\ 0 & 0 \end{pmatrix} \begin{pmatrix} M_{11} & M_{12} \\ M_{21} & M_{22} \end{pmatrix}^T = \begin{pmatrix} M_{11}M_{11}^T & M_{11}M_{21}^T \\ M_{21}M_{11}^T & M_{21}M_{21}^T \end{pmatrix}.
\]
It follows that $M_{21} = 0$ and $M_{11}$ is an $r \times r$ orthogonal matrix. However, it is well-known that orthogonal matrices possess a full basis of eigenvectors with eigenvalues on the unit circle. The assertion of the lemma now readily follows.
\end{proof}

%
%

\section{Copositive matrices and linear periodic systems}

In this section we establish a relation between the objects considered in the preceding two sections. Throughout this and the next section, we fix a collection ${\bf u} = \{u^1,\dots,u^n\} \subset \mathbb R_+^n$ of nonnegative vectors such that $\Supp u^j = I_j$, $j = 1,\dots,n$. Moreover, we assume these vectors are normalized such that the last elements of their positive subvectors $u^j_{I_j}$ all equal 1. With the collection $\bf u$ we associate a discrete-time linear periodic system ${\bf S}_{\bf u}$ of order $d = n-3$ and with period $n$, given by \eqref{periodic_system} with coefficient vectors $c^t = u^t_{I_t}$, $t = 1,\dots,n$. The coefficient vectors $c^t$ for all other time instants $t > n$ are then determined by the periodicity relation $c^{t+n} = c^t$. Equivalently, the dynamics of ${\bf S}_{\bf u}$ is given by the equations
\begin{equation} \label{def_system_Su}
\sum_{i=0}^d u^{t'}_sx_{t+i} = 0,\qquad t \geq 1,
\end{equation}
where $t',s \in \{1,\dots,n\}$ are the unique indices such that $t \equiv t'$ and $t+i \equiv s$ modulo $n$. Relation \eqref{e_relation} then becomes
\begin{equation} \label{e_relation2}
\sum_{i=0}^d u^{t'}_s{\bf e}_{t+i} = 0,\qquad \forall\ t \geq 1,
\end{equation}
with $t',s$ defined as above.

By Lemma \ref{monodromy_det} the monodromy of ${\bf S}_{\bf u}$ then satisfies
\begin{equation} \label{det_positive}
\det {\mathfrak M} = \prod_{j=1}^n u^j_j > 0.
\end{equation}
In particular, the system ${\bf S}_{\bf u}$ is time-reversible. Denote by ${\cal L}_{\bf u}$ the space of solutions of ${\bf S}_{\bf u}$.

Let ${\cal A}_{\bf u} \subset {\cal S}^n$ be the linear subspace of matrices $A$ satisfying $A_{I_j}u^j_{I_j} = A_{I_j}c^j = 0$ for all $j = 1,\dots,n$. With $A \in {\cal A}_{\bf u}$ we associate a symmetric bilinear form $B$ on the dual space ${\cal L}_{\bf u}^*$ by setting $B({\bf e}_t,{\bf e}_s) = A_{ts}$ for every $t,s = 1,\dots,d$ and defining the value of $B$ on arbitrary vectors in ${\cal L}_{\bf u}^*$ by linear extension. In other words, in the basis $\{{\bf e}_1,\dots,{\bf e}_d\}$ of ${\cal L}_{\bf u}^*$ the coefficient matrix of $B$ is given by the submatrix $A_{I_1'}$. Let $\Lambda: A \mapsto B$ be the so-defined linear map from ${\cal A}_{\bf u}$ into the space of symmetric bilinear forms on ${\cal L}_{\bf u}^*$. Our first step will be to describe the image of $\Lambda$. To this end, we need the following lemma.

{\lemma \label{lem:imageL} Let $A \in {\cal A}_{\bf u}$ and $B = \Lambda(A)$. Then for every integer $r \geq 1$, the $(n-2) \times (n-2)$-matrix $B_r = (B({\bf e}_{t+r},{\bf e}_{s+r}))_{t,s = 0,\dots,d}$ equals the submatrix $A_{I_{r'}}$, where $r' \in \{1,\dots,n\}$ is the unique index satisfying $r \equiv r'$ modulo $n$. Equivalently,
\begin{equation} \label{band_relation}
B({\bf e}_t,{\bf e}_s) = A_{t's'}\qquad \forall\ t,s \geq 1:\ |t - s| \leq n-3,
\end{equation}
where $t',s' \in \{1,\dots,n\}$ are the unique indices such that $t \equiv t'$, $s'\equiv s$ modulo $n$. }

\begin{proof}
We proceed by induction over $r$. By definition of $\Lambda$ the upper left $(n-3) \times (n-3)$ submatrix of $B_1$ equals the corresponding submatrix of $A_{I_1}$. However, we have $A_{I_1}c^1 = 0$ by virtue of $A \in {\cal A}_{\bf u}$. Moreover, $B_1c^1 = 0$, which can be seen as follows. For arbitrary $i = 0,\dots,d$ we have
\[ (B_1c^1)_{i+1} = \sum_{j=0}^dB({\bf e}_{i+1},{\bf e}_{j+1})c^1_j = B({\bf e}_{i+1},\sum_{j=0}^d{\bf e}_{j+1}c^1_j) = B({\bf e}_{i+1},0) = 0
\]
by virtue of \eqref{e_relation} for $t = 1$. Hence the difference $A_{I_1} - B_1$ is a symmetric matrix, with possibly non-zero elements only in the last row or in the last column, which possesses a kernel vector with all elements positive. It is easily seen that this difference must then be the zero matrix, proving the assertion of the lemma for $r = 1$.

The induction step from $r-1$ to $r$ proceeds in a similar manner, with equality of the upper left $(n-3) \times (n-3)$ submatrices of $B_r$ and $A_{I_{r'}}$ now guaranteed by the induction hypothesis.
\end{proof}

The lemma asserts that the submatrices $A_{I_1},\dots,A_{I_n}$ are all of the form $(B({\bf e}_t,{\bf e}_s))_{t,s \in I}$ for certain index sets $I$. Some of their properties are hence determined by the corresponding properties of $B$.

{\corollary \label{cor:rank_equal} Let $A \in {\cal A}_{\bf u}$ and $B = \Lambda(A)$. Then the ranks of the matrices $A_{I_j}$, $j = 1,\dots,n$, and of all their submatrices of size $(n-3) \times (n-3)$, are equal to the rank of the symmetric bilinear form $B$. }

\begin{proof}
Since the system ${\bf S}_{\bf u}$ is time-reversible, the evaluation operators ${\bf e}_t,\dots,{\bf e}_{t+d-1}$ form a basis of ${\cal L}_{\bf u}^*$ for every $t \geq 1$. On the other hand, the operators ${\bf e}_t,\dots,{\bf e}_{t+d}$ are linearly dependent for all $t \geq 1$, the dependence being given by \eqref{e_relation2}. All coefficients in this relation are non-zero, hence any of the operators ${\bf e}_t,\dots,{\bf e}_{t+d}$ can be expressed as a linear combination of the $d$ other operators. It follows that every subset of $\{{\bf e}_t,\dots,{\bf e}_{t+d}\}$ of cardinality $d$ is a basis of ${\cal L}_{\bf u}^*$. Therefore the $d \times d$ matrix $(B({\bf e}_s,{\bf e}_{s'}))_{s \in I,s'\in I'}$, where $I,I'$ are such subsets, has rank equal to $\rk\,B$. The same holds for the matrices $(B({\bf e}_s,{\bf e}_{s'}))_{s,s' = t,\dots,t+d}$ for all $t \geq 1$. The corollary now follows from Lemma \ref{lem:imageL}.
\end{proof}

{\corollary \label{cor:imageL} Let $A \in {\cal A}_{\bf u}$ and $B = \Lambda(A)$. Then the symmetric bilinear form $B$ is shift-invariant and satisfies the linear relations
\begin{equation} \label{linear_relation}
B({\bf e}_t,{\bf e}_s) = B({\bf e}_{t+n},{\bf e}_s)\qquad \forall\ t,s \geq 1:\ 3 \leq s-t \leq n-3.
\end{equation}}

\begin{proof}
By \eqref{band_relation} we have $B({\bf e}_t,{\bf e}_s) = A_{ts} = B({\bf e}_{t+n},{\bf e}_{s+n})$ for all $t,s = 1,\dots,d$. This in turn implies the shift-invariance of $B$ by Lemma \ref{lem:inv}.

The inequalities $3 \leq s-t \leq n-3$ imply $|t - s| \leq n-3$, $|t + n - s| \leq n-3$. Relations \eqref{linear_relation} then follow from \eqref{band_relation} in a similar way as the shift-invariance.
\end{proof}

Now we are ready to describe the image of the map $\Lambda$.

{\lemma \label{lem:injective} Suppose that $n \geq 5$. Then the linear map $\Lambda$ is injective, and its image consists of those shift-invariant symmetric bilinear forms $B$ on ${\cal L}_{\bf u}^*$ which satisfy relations \eqref{linear_relation}. }

\begin{proof}
In view of Corollary \ref{cor:imageL} it suffices to show that for every shift-invariant symmetric bilinear form $B$ satisfying \eqref{linear_relation} there exists a unique matrix $A \in {\cal A}_{\bf u}$ such that $B = \Lambda(A)$.

Let $B$ be such a form. We define the corresponding matrix $A$ as follows. Let $i,j \in \{1,\dots,n\}$ be arbitrary indices such that $i \leq j$. Put
\begin{equation} \label{Lambda_explicit}
A_{ij} = A_{ji} = \left\{ \begin{array}{rcl} B({\bf e}_i,{\bf e}_j),& \quad & j - i \leq n-3; \\ B({\bf e}_{i+n},{\bf e}_j), && j - i > n-3. \end{array} \right.
\end{equation}
The shift-invariance of $B$ and relations \eqref{linear_relation} then imply \eqref{band_relation}. This yields $B_r = (B({\bf e}_{t+r},{\bf e}_{s+r}))_{t,s = 0,\dots,d} = A_{I_r}$ for every $r = 1,\dots,n$. Now $B_rc^r = 0$ by virtue of \eqref{e_relation}, which implies $A_{I_r}c^r = 0$ and hence $A \in {\cal A}_{\bf u}$. Moreover, $\Lambda(A) = B$ by construction of $A$.

Uniqueness of $A$ follows from Lemmas \ref{lem:imageL} and \ref{lem:irred_N}.
\end{proof}

\medskip

Now we shall investigate which symmetric bilinear forms $B$ in the image of $\Lambda$ are the images of copositive matrices. First we consider positive semi-definite bilinear symmetric forms $B$.

{\lemma \label{lem:Bpos} Suppose $n \geq 5$ and let $A \in {\cal A}_{\bf u}$ and $B = \Lambda(A)$. Then the following are equivalent:
\begin{itemize}
\setlength{\itemsep}{1pt}
\setlength{\parskip}{0pt}
\setlength{\parsep}{0pt}
\item[(i)] the form $B$ is positive semi-definite;
\item[(ii)] the submatrices $A_{I_j}$ are positive semi-definite for all $j = 1,\dots,n$;
\item[(iii)] any of the submatrices $A_{I_j}$, $j = 1,\dots,n$, is positive semi-definite.
\end{itemize}

Moreover, given above conditions (i)---(iii), the following holds:
\begin{itemize}
\setlength{\itemsep}{1pt}
\setlength{\parskip}{0pt}
\setlength{\parsep}{0pt}
\item[(a)] the difference $B({\bf e}_n,{\bf e}_{n-2}) - B({\bf e}_n,{\bf e}_{2n-2})$ has the same sign as $(u^n)^TAu^1$;
\item[(b)] the difference $B({\bf e}_{j+n},{\bf e}_{j+n-2}) - B({\bf e}_{j+n},{\bf e}_{j+2n-2})$ has the same sign as $(u^j)^TAu^{j+1}$, $j = 1,2$;
\item[(c)] the difference $B({\bf e}_j,{\bf e}_{j-2}) - B({\bf e}_j,{\bf e}_{j+n-2})$ has the same sign as $(u^j)^TAu^{j+1}$, $j = 3,\dots,n-1$.
\end{itemize} }

\begin{proof}
The first part of the lemma is a direct consequence of Lemma \ref{lem:imageL} and of the fact that the set $\{ {\bf e}_t,\dots,{\bf e}_{t+d} \}$ spans the whole space ${\cal L}_{\bf u}^*$ for all $t \geq 1$.

In order to prove the second part, let us assume conditions (i)---(iii). Consider the $(n-1) \times (n-1)$ matrix $A_{(1,\dots,n-1)}$. By Lemma \ref{lem:imageL} its upper left and its lower right principal submatrix of size $n-2$ coincides with the matrix $(B({\bf e}_t,{\bf e}_s))_{t,s \in I}$ with $I = \{1,\dots,n-2\}$ and $I = \{2,\dots,n-1\}$, respectively. Hence it can be written as a sum $(B({\bf e}_t,{\bf e}_s))_{t,s = 1,\dots,n-1} + \delta E_{1,n-1}$ for some real $\delta$. Note that the first summand is positive semi-definite by condition (i). On the other hand, $A_{1,n-1} = B({\bf e}_{n+1},{\bf e}_{n-1})$ by \eqref{band_relation}. Hence $\delta = B({\bf e}_{n+1},{\bf e}_{n-1}) - B({\bf e}_1,{\bf e}_{n-1}) = B({\bf e}_{n+1},{\bf e}_{n-1}) - B({\bf e}_{n+1},{\bf e}_{2n-1})$, where the second equality follows from the shift-invariance of $B$ which is in turn a consequence of Corollary \ref{cor:imageL}. By virtue of Lemma \ref{completion} we then get (b) for $j = 1$.

The other assertions of the second part are proven in a similar way by starting with the remaining $(n-1) \times (n-1)$ principal submatrices of $A$.
\end{proof}

%
%
%

{\lemma \label{lem:second_cond} Let $n \geq 5$, and let $A \in {\cal C}^n$ be such that $(u^j)^TAu^j = 0$ for all $j = 1,\dots,n$. Then $A \in {\cal A}_{\bf u}$, and $B = \Lambda(A)$ is positive semi-definite and satisfies the inequalities
\begin{equation} \label{ineq_relations}
B({\bf e}_t,{\bf e}_{t+2}) \geq B({\bf e}_{t+n},{\bf e}_{t+2})\qquad \forall t \geq 1.
\end{equation}
Moreover, either $B$ satisfies all inequalities \eqref{ineq_relations} with equality, namely when $A$ is positive semi-definite, or all inequalities \eqref{ineq_relations} are strict, namely when $A$ is exceptional. }

\begin{proof}
Let $A \in {\cal C}^n$ such that $(u^j)^TAu^j = 0$ for all $j = 1,\dots,n$. Then $A_{I_j} \succeq 0$ for $j = 1,\dots,n$ by Lemma \ref{zeroPSD}. The relation $(c^j)^TA_{I_j}c^j = (u^j)^TAu^j = 0$ then implies $A_{I_j}c^j = 0$ for all $j = 1,\dots,n$, and we get indeed $A \in {\cal A}_{\bf u}$. By Lemma \ref{lem:Bpos} we then have $B \succeq 0$. From Theorem \ref{first_cond} and Lemma \ref{lem:Bpos} we obtain \eqref{ineq_relations} for $t = 1,\dots,n$, with equality or strict inequality for positive semi-definite or exceptional $A$, respectively. For all other $t > n$ these relations follow by the shift-invariance of $B$ which in turn is implied by Corollary \ref{cor:imageL}.
\end{proof}

{\lemma \label{lem:other_side} Let $B$ be a positive semi-definite symmetric bilinear form in the image of $\Lambda$ which satisfies inequalities \eqref{ineq_relations}. Then its pre-image $A = \Lambda^{-1}(B)$ is copositive and satisfies $(u^j)^TAu^j = 0$ for all $j = 1,\dots,n$. }

\begin{proof}
Let $B$ be as required, and let $A$ be its pre-image, which is well-defined by Lemma \ref{lem:injective}. By Lemma \ref{lem:Bpos} the submatrices $A_{I_j}$ are positive semi-definite for all $j = 1,\dots,n$. From $A \in {\cal A}_{\bf u}$ it follows that $(u^j)^TAu^j = 0$ for all $j = 1,\dots,n$. Finally, by Lemma \ref{lem:Bpos} \eqref{ineq_relations} implies the inequalities $(u^n)^TAu^1 \geq 0$ and $(u^j)^TAu^{j+1} \geq 0$, $j = 1,\dots,n-1$. Therefore $A \in {\cal C}^n$ by Theorem \ref{first_cond}.
\end{proof}

Now we are in a position to describe the face $F_{\bf u}$ of the copositive cone and the face $P_{\bf u}$ of the positive semi-definite cone which are defined by virtue of \eqref{faces_def} by the zeros $u^j$, $j = 1,\dots,n$. The description will be by linear matrix inequalities.

{\theorem \label{thm:sdp_rep} Let $n \geq 5$, and let ${\cal F}_{\bf u}$ be the set of positive semi-definite symmetric bilinear forms $B$ on ${\cal L}_{\bf u}^*$ satisfying the linear equality relations
\[
\begin{array}{rclcl}
B({\bf e}_t,{\bf e}_s) &=& B({\mathfrak M}^*{\bf e}_t,{\mathfrak M}^*{\bf e}_s),&\qquad& t,s = 1,\dots,n-3; \\
B({\bf e}_t,{\bf e}_s) &=& B({\mathfrak M}^*{\bf e}_t,{\bf e}_s),&\qquad& 1 \leq t < s \leq n:\ 3 \leq s-t \leq n-3
\end{array}
\]
and the linear inequalities
\[ B({\bf e}_t,{\bf e}_{t+2}) \geq B({\mathfrak M}^*{\bf e}_t,{\bf e}_{t+2}),\qquad t = 1,\dots,n.
\]
Let ${\cal P}_{\bf u} \subset {\cal F}_{\bf u}$ be the subset of forms $B$ which satisfy all linear inequalities with equality.

Then the face of ${\cal C}^n$ defined by the zeros $u^j$, $j = 1,\dots,n$, satisfies $F_{\bf u} = \Lambda^{-1}[{\cal F}_{\bf u}]$, and the face of ${\cal S}_+^n$ defined by these zeros satisfies $P_{\bf u} = \Lambda^{-1}[{\cal P}_{\bf u}]$. Moreover, for all forms $B \in {\cal F}_{\bf u} \setminus {\cal P}_{\bf u}$ the linear inequalities are satisfied strictly, and $F_{\bf u} \setminus P_{\bf u}$ consists of exceptional matrices. }


\begin{proof}
By virtue of Lemmas \ref{lem:injective}, \ref{lem:second_cond}, and \ref{lem:other_side} we have only to show that the finite number of equalities and inequalities stated in the formulation of the theorem are necessary and sufficient to ensure the infinite number of equalities and inequalities in \eqref{periodic_Q},\eqref{linear_relation},\eqref{ineq_relations}. Necessity is evident, since the relations in the theorem are a subset of relations \eqref{periodic_Q},\eqref{linear_relation},\eqref{ineq_relations}.

Sufficiency of the first set of equalities in the theorem to ensure \eqref{periodic_Q} follows from Lemma \ref{lem:inv}. By the resulting shift-invariance of $B$ it is also sufficient to constrain the index $t$ to $1,\dots,n$ in \eqref{linear_relation}. Now suppose that $t \in \{1,\dots,n\}$, $s > n$, and $3 \leq s-t \leq n-3$. We then have $B({\bf e}_t,{\bf e}_s) = B({\bf e}_s,{\bf e}_t)$, $B({\bf e}_{t+n},{\bf e}_s) = B({\bf e}_t,{\bf e}_{s-n}) = B({\bf e}_{s-n},{\bf e}_t)$. Relation \eqref{linear_relation} for the index pair $(t,s)$ is hence equivalent to the same relation for the index pair $(t',s') = (s-n,t)$, which also satisfies $1 \leq t' \leq n$ and $3 \leq s'-t' \leq n-3$. For the latter index pair we now have $s' \leq n$, however. Thus \eqref{linear_relation} also follows from the linear equalities in the theorem. Finally, the set of inequalities in the theorem implies \eqref{ineq_relations} by the shift-invariance of $B$.
\end{proof}

{\corollary \label{cor:rank_drop} The rank of the forms in ${\cal F}_{\bf u}$ is constant over the (relative) interior of ${\cal F}_{\bf u}$ and there it is equal to the maximum of the rank over all forms in ${\cal F}_{\bf u}$. Let us the denote this maximal rank by $r_{\max}$. Then for every $B \in \partial {\cal F}_{\bf u} \setminus {\cal P}_{\bf u}$ we have $\rk\,B < r_{\max}$. }

\begin{proof}
By definition ${\cal F}_{\bf u}$ is the intersection of the cone of positive semi-definite symmetric bilinear forms on ${\cal L}_{\bf u}^*$ with a homogeneous polyhedral set. The relative interior of ${\cal F}_{\bf u}$ is hence non-empty and contained in the relative interior of some face of this cone. Therefore the rank is constant over the relative interior of ${\cal F}_{\bf u}$. Moreover, since the rank is a lower semi-continuous function, it can only drop on the relative boundary of ${\cal F}_{\bf u}$. This proves the first claim.

Let now $B \in \partial {\cal F}_{\bf u} \setminus {\cal P}_{\bf u}$, and let $B'$ be a form in the relative interior of ${\cal F}_{\bf u}$. Let $l$ be the line segment connecting $B$ with $B'$. When we reach the boundary point $B$ by moving along $l$, either a rank drop must occur, in case the relative boundary of the face of the positive semi-definite cone is reached, or one of the linear inequalities must become an equality. However, the second case cannot happen, because $B \not\in {\cal P}_{\bf u}$.
\end{proof}

If the vectors in the collection $\bf u$ are in generic position, then the set ${\cal F}_{\bf u}$ defined in Theorem \ref{thm:sdp_rep} consists of the zero form only. In the next section we investigate the consequences of a non-trivial set ${\cal F}_{\bf u}$.

\section{Structure of the cones ${\cal F}_{\bf u}$ and ${\cal P}_{\bf u}$} \label{sec:structure}

As was mentioned in the introduction, the eigenvalues of the monodromy $\mathfrak M$, the Floquet multipliers, largely determine the properties of the matrices in the face $F_{\bf u}$ of ${\cal C}^n$. In this section we shall investigate these connections in detail. In particular, we will be interested in the structure of the cones $P_{\bf u}$ and ${\cal P}_{\bf u} = \Lambda[P_{\bf u}]$ defined by the positive semi-definite matrices in the face $F_{\bf u} \subset {\cal C}^n$ and its connections to the periodic solutions of the system ${\bf S}_{\bf u}$. We also investigate the properties of the exceptional copositive matrices with minimal and non-minimal circulant zero support set as defined in Definition \ref{def:ref}.


Denote by ${\cal L}_{per} \subset {\cal L}_{\bf u}$ the subspace of $n$-periodic solutions. We have the following characterization of ${\cal L}_{per}$.

{\lemma \label{periodic_char} An $n$-periodic infinite sequence $x = (x_1,x_2,\dots)$ is a solution of ${\bf S}_{\bf u}$ if and only if the vector $(x_1,\dots,x_n)^T \in \mathbb R^n$ is orthogonal to all vectors $u^j$, $j = 1,\dots,n$. In particular, the dimension of ${\cal L}_{per}$ equals the corank of the $n \times n$ matrix $U$ composed of the column vectors $u^1,\dots,u^n$. }

\begin{proof}
From the $n$-periodicity of $x$ it follows that \eqref{def_system_Su} are exactly the orthogonality relations between $(x_1,\dots,x_n)^T$ and $u^j$. The lemma now readily follows.
\end{proof}

We are now in a position to describe the set ${\cal P}_{\bf u}$ in terms of the subspace ${\cal L}_{per}$.

{\lemma \label{eigenvalue1} Suppose that $n \geq 5$. Then ${\cal P}_{\bf u}$ equals the convex hull of all tensor products $x \otimes x$, $x \in {\cal L}_{per}$. In particular, ${\cal P}_{\bf u} \simeq {\cal S}_+^{\dim{\cal L}_{per}}$, and for every $B \in {\cal P}_{\bf u}$ we have $\it{Im}\,B \subset {\cal L}_{per}$. Moreover, for every $B \in {\cal P}_{\bf u}$ the preimage $A = \Lambda^{-1}(B)$ is given by $A = (B({\bf e}_t,{\bf e}_s))_{t,s = 1,\dots,n}$. }

\begin{proof}
Assume $x \in {\cal L}_{per}$ and set $B = x \otimes x$. The $n$-periodicity of the solution $x$ implies that $x_{t+n} = x_t$ for all $t \geq 1$. For every $t,s \geq 1$ we then have $B({\bf e}_t,{\bf e}_s) = x_tx_s = x_{t+n}x_s = B({\bf e}_{t+n},{\bf e}_s)$, which yields \eqref{linear_relation} and \eqref{ineq_relations} with equality. In a similar way we obtain \eqref{periodic_Q}, and hence $B \in {\cal P}_{\bf u}$. Moreover, \eqref{Lambda_explicit} yields $A = \Lambda^{-1}(B) = (B({\bf e}_t,{\bf e}_s))_{t,s = 1,\dots,n}$.

By convexity of ${\cal P}_{\bf u}$ it follows that the convex hull of the set $\{ x \otimes x \,|\, x \in {\cal L}_{per} \}$ is a subset of ${\cal P}_{\bf u}$, and by linearity the above expression for $A = \Lambda^{-1}(B)$ holds also for every $B$ in this convex hull.

Let us prove the converse inclusion. Let $B \in {\cal P}_{\bf u}$ be arbitrary and set $A = \Lambda^{-1}(B)$. Then $A \in P_{\bf u}$ by Theorem \ref{thm:sdp_rep}. Since $A \succeq 0$ and $(u^j)^TAu^j = 0$, it follows that $Au^j = 0$ for all $j = 1,\dots,n$. Therefore $A$ is in the convex hull of the set $\{ vv^T \,|\, \langle v,u^j \rangle = 0\ \forall\ j = 1,\dots,n \}$. It hence suffices to show that for every $v \in \mathbb R^n$ such that $\langle v,u^j \rangle = 0$, $j = 1,\dots,n$, we have $\Lambda(vv^T) = x \otimes x$ for some $x \in {\cal L}_{per}$.

Let $v \in \mathbb R^n$ be orthogonal to all zeros $u^j$. By Lemma \ref{periodic_char} the $n$-periodic infinite sequence  $x = (x_1,x_2,\dots)$ defined by the initial conditions $x_i = v_i$, $i = 1,\dots,n$, is an $n$-periodic solution of ${\bf S}_{\bf u}$, $x \in {\cal L}_{per}$. But $\Lambda(vv^T) = x \otimes x$ by construction. This completes the proof.
\end{proof}

We now consider the ranks of the forms in ${\cal F}_{\bf u}$ and ${\cal P}_{\bf u}$. Let $r_{\max}$ be the maximal rank achieved by forms in ${\cal F}_{\bf u}$, and $r_{\it{PSD}}$ the maximal rank achieved by forms in ${\cal P}_{\bf u}$.

{\corollary \label{cor:eigenvalue1} Suppose $n \geq 5$. Then $r_{\it{PSD}}$ equals the geometric multiplicity of the eigenvalue 1 of the monodromy operator $\mathfrak M$ of the dynamical system ${\bf S}_{\bf u}$. }

\begin{proof}
By Lemma \ref{eigenvalue1} we have $r_{PSD} = \dim {\cal L}_{per}$. The corollary then follows from Lemma \ref{periodic1}.
\end{proof}

{\lemma \label{lem:triangular} Let $n \geq 5$, and let $B \in {\cal F}_{\bf u} \setminus {\cal P}_{\bf u}$. Then for every $k \geq 1$ the matrix
\[ M_k = (B({\bf e}_{t+k}-{\bf e}_{t+n+k},{\bf e}_{s+k}))_{t = 0,\dots,n-5;s = 2,\dots,n-3} = (B((\it{Id} - {\mathfrak M}^*){\bf e}_{t+k},{\bf e}_{s+k}))_{t = 0,\dots,n-5;s = 2,\dots,n-3}
\]
has full rank $n-4$. }

\begin{proof}
By Lemma \ref{lem:second_cond} the form $B$ satisfies inequalities \eqref{ineq_relations} strictly, which implies that $M_k$ has all diagonal elements positive. By \eqref{linear_relation} $M_k$ is lower-triangular, and hence has full rank $n-4$.
\end{proof}

%

{\corollary \label{cor:corank1} Let $n \geq 5$, and let $B \in {\cal F}_{\bf u} \setminus {\cal P}_{\bf u}$. Then the bilinear form on ${\cal L}_{\bf u}^*$ given by $(w,w') \mapsto B((\it{Id}-{\mathfrak M}^*)w,w')$ has corank at most 1. In particular, both $B$ and $\it{Id}-{\mathfrak M}^*$ have corank at most 1. Moreover, $\mathfrak M$ has at least $n-4$ linearly independent eigenvectors with eigenvalues on the unit circle. }

\begin{proof}
The bilinear form in the statement of the lemma has at least the same rank, namely $n-4$, as the matrices $M_k$ in Lemma \ref{lem:triangular}. Hence it has corank at most 1. It follows that $B$ has corank at most 1, and the proof is concluded by application of Lemma \ref{lem:inv_unitary}.
\end{proof}

{\corollary \label{cor:difference_ranks} Suppose $n \geq 5$. If ${\cal F}_{\bf u} \not= {\cal P}_{\bf u}$, then $r_{\max} - r_{\it{PSD}} \geq n-4$. In particular, in this case $r_{\it{PSD}} \leq 1$ and either $r_{\max} = n - 4$ or $r_{\max} = n - 3$. }

\begin{proof}
Let $B \in {\cal F}_{\bf u} \setminus {\cal P}_{\bf u}$. Suppose there exists $B' \in {\cal P}_{\bf u}$ such that $B \succeq B'$. Then also $B - B' \in {\cal F}_{\bf u} \setminus {\cal P}_{\bf u}$. Therefore we may assume without loss of generality that there does not exist a non-zero $B' \in {\cal P}_{\bf u}$ such that $B - B' \succeq 0$. By Lemma \ref{eigenvalue1} we then have $\it{Im}\,B \cap {\cal L}_{per} = \{0\}$, and hence for every $B' \in {\cal P}_{\bf u}$ we get $\rk(B+B') = \rk\,B + \rk\,B'$, again by Lemma \ref{eigenvalue1}.

Let now $B' \in {\cal P}_{\bf u}$ such that $\rk\,B' = r_{\it{PSD}}$. By Corollary \ref{cor:corank1} we have $\rk\,B \geq n-4$, and hence $r_{\max} \geq \rk\,B + \rk\,B' \geq n - 4 + r_{\it{PSD}}$. This completes the proof.
\end{proof}

%
%

These results allow us to completely characterize the face $F_{\bf u}$ in the case when ${\cal F}_{\bf u}$ does not contain positive definite forms.

{\lemma \label{exc_degenerate} Suppose $n \geq 5$ and assume $r_{\max} = n - 4$. Then either $F_{\bf u}$ consists of positive semi-definite matrices only, or $F_{\bf u}$ is 1-dimensional and generated by an extremal exceptional copositive matrix $A$. In the latter case the submatrices $A_{I_j}$ of this exceptional matrix have corank 2 for all $j = 1,\dots,n$. }

\begin{proof}
By Corollary \ref{cor:rank_drop} we have $\rk B \leq n-5$ for every $B \in \partial {\cal F}_{\bf u} \setminus {\cal P}_{\bf u}$. By Corollary \ref{cor:corank1}, however, the corank of any such matrix can be at most 1. Therefore $\partial {\cal F}_{\bf u} \setminus {\cal P}_{\bf u} = \emptyset$, and we have the inclusion $\partial {\cal F}_{\bf u} \subset {\cal P}_{\bf u}$. But ${\cal P}_{\bf u}$ is a cone, and hence it contains also the conic hull of $\partial {\cal F}_{\bf u}$.

However, the conic hull of the boundary of some pointed cone is either the cone itself, namely when the dimension of the cone exceeds 1, or it is the singleton $\{0\}$, namely when the cone is 1-dimensional. Therefore if ${\cal F}_{\bf u} \not= {\cal P}_{\bf u}$, then ${\cal F}_{\bf u}$ must be 1-dimensional and we must have ${\cal P}_{\bf u} = \{0\}$.

The first claim of the lemma now follows from Theorem \ref{thm:sdp_rep}. The second claim then follows from Corollary \ref{cor:rank_equal}.
\end{proof}

We shall now concentrate on the case when ${\cal F}_{\bf u}$ contains positive definite forms, i.e., the maximal rank achieved by matrices in ${\cal F}_{\bf u}$ equals $r_{\max} = d = n-3$.

{\lemma \label{PSD_full_rank} Let $n \geq 5$. Then the following are equivalent:
\begin{itemize}
\setlength{\itemsep}{1pt}
\setlength{\parskip}{0pt}
\setlength{\parsep}{0pt}
\item[(i)] the set ${\cal F}_{\bf u}$ contains a positive definite form and ${\cal F}_{\bf u} = {\cal P}_{\bf u}$;
\item[(ii)] the monodromy operator $\mathfrak M$ of the system ${\bf S}_{\bf u}$ equals the identity;
\item[(iii)] the rank of the $n \times n$ matrix $U$ with columns $u^1,\dots,u^n$ equals 3.
\end{itemize} }

\begin{proof}
Assume (ii). The LMIs in Theorem \ref{thm:sdp_rep} reduce to the condition $B \succeq 0$, and hence ${\cal F}_{\bf u} = {\cal P}_{\bf u} \simeq {\cal S}_+^{n-3}$, implying (i).

Assume (i). We have $r_{\max} = r_{\it{PSD}} = n-3$, and (ii) follows from Corollary \ref{cor:eigenvalue1}.

By Lemma \ref{periodic1} we have ${\mathfrak M} = \it{Id}$ if and only if $\dim{\cal L}_{per} = n-3$. By Lemma \ref{periodic_char} this is equivalent to the condition $\rk U = 3$, which proves (ii) $\Leftrightarrow$ (iii).
\end{proof}

In order to treat the case $r_{\max} = n - 3$ and ${\cal F}_{\bf u} \not= {\cal P}_{\bf u}$ we shall distinguish between odd and even orders $n$.

{\lemma \label{full_rank_even} Let $n > 5$ be even, and suppose $r_{\max} = n - 3$ and ${\cal F}_{\bf u} \not= {\cal P}_{\bf u}$. Then $F_{\bf u}$ is linearly isomorphic to $\mathbb R_+^2$, where one boundary ray of $F_{\bf u}$ is generated by a rank 1 positive semi-definite matrix, and the other boundary ray is generated by an extremal exceptional copositive matrix $A$. The submatrices $A_{I_j}$ of this exceptional matrix have corank 2 for all $j = 1,\dots,n$. }

\medskip

Before we proceed to the proof, we shall provide the following fact from linear algebra.

\emph{Let $V$ be a real vector space, let $A: V \to V$ be a linear map, and let $A^*: V^* \to V^*$ be its adjoint. Let $\lambda \in \mathbb R$ be an eigenvalue of $A$ and $A^*$, and let $U \subset V$, $W \subset V^*$ be the corresponding invariant subspaces of $A$ and $A^*$, respectively. Then $V^*$ can be written as a direct sum $W \oplus U^{\perp}$, where $U^{\perp}$ is the annihilator of $U$. }

The claim becomes immediate if we pass to a coordinate system where $A$ is in real canonical Jordan form.

\begin{proof}(of Lemma \ref{full_rank_even})
By Lemma \ref{lem:inv_unitary} all eigenvalues of the monodromy $\mathfrak M$ of the system ${\bf S}_{\bf u}$ lie on the unit circle and their geometric and algebraic multiplicities coincide. However, since $\mathfrak M$ is real, its complex eigenvalues are grouped into complex-conjugate pairs. Since the dimension $d = n-3$ of ${\cal L}_{\bf u}^*$ is odd, there must be exactly one real eigenvalue with an odd multiplicity. By \eqref{det_positive} this eigenvalue equals 1. Hence $r_{\it{PSD}}$ is odd by Corollary \ref{cor:eigenvalue1}, but cannot exceed 1 by Corollary \ref{cor:difference_ranks}. Therefore by Corollary \ref{cor:eigenvalue1} and Lemma \ref{periodic1} we have $\dim{\cal L}_{per} = 1$. By Lemma \ref{eigenvalue1} we have $\dim {\cal P}_{\bf u} = \dim P_{\bf u} = 1$, and $P_{\bf u}$ is generated by a rank 1 positive semi-definite matrix $A_P$.

Denote the 1-dimensional eigenspace of $\mathfrak M^*$ corresponding to the eigenvalue 1 by $W_1$, and let ${\cal L}_{per}^{\perp} \subset {\cal L}_{\bf u}^*$ be the annihilator of ${\cal L}_{per}$. Then ${\cal L}_{per}^{\perp}$ is an invariant subspace of $\mathfrak M^*$ and by the above ${\cal L}_{\bf u}^* = W_1 \oplus {\cal L}_{per}^{\perp}$. Let now $w^1 \in W_1$ and $w \in {\cal L}_{per}^{\perp}$ be arbitrary vectors. Then for every $B \in {\cal F}_{\bf u}$ we get by the shift-invariance $B(w^1,w) = B(\mathfrak M^*w^1,\mathfrak M^*w) = B(w^1,\mathfrak M^*w)$, and hence $B(w^1,(\it{Id} - \mathfrak M^*)w) = 0$ for all $w \in {\cal L}_{per}^{\perp}$. But $(\it{Id} - \mathfrak M^*)[{\cal L}_{per}^{\perp}] = {\cal L}_{per}^{\perp}$, because ${\cal L}_{per}^{\perp}$ is an invariant subspace of $\it{Id} - \mathfrak M^*$ and it has a zero intersection with the kernel\footnote{The kernel of a symmetric bilinear form $B$ is the set of vectors $x$ such that $B(x,y) = 0$ for all $y$.} $W_1$ of $\it{Id} - \mathfrak M^*$. It follows that $W_1$ and ${\cal L}_{per}^{\perp}$ are orthogonal under $B$.

For every pair of vectors $w^1 \in W_1$, $w \in {\cal L}_{per}^{\perp}$ we then get $B(w^1+w,w^1+w) = B(w^1,w^1) + B(w,w)$. This implies that every $B \in {\cal F}_{\bf u}$ can in a unique way be decomposed into a sum $B = B' + P$ of positive semi-definite bilinear symmetric forms, with $P \in {\cal P}_{\bf u}$ and $W_1 \subset \ker B'$. Namely, $B',P$ are defined such that $P(v,v) = B(w^1,w^1)$, $B'(v,v) = B(w,w)$ for every $v \in {\cal L}_{\bf u}^*$, where $v = w^1 + w$ is the unique decomposition of $v$ into vectors $w^1 \in W_1$, $w \in {\cal L}_{per}^{\perp}$. Moreover, for every $B \in {\cal F}_{\bf u}$ the corresponding summand $B'$ is also in ${\cal F}_{\bf u}$, because $P$ satisfies inequalities \eqref{ineq_relations} with equality. Thus the cone ${\cal F}_{\bf u}$ splits into a direct sum ${\cal F}_{\bf u}' + {\cal P}_{\bf u}$, where ${\cal F}_{\bf u}' = \{ B \in {\cal F}_{\bf u} \,|\, W_1 \subset \ker B \}$.

By assumption ${\cal F}_{\bf u}' \not= \{0\}$. Any non-zero form in ${\cal F}_{\bf u}'$ lies in $\partial {\cal F}_{\bf u} \setminus {\cal P}_{\bf u}$ and hence must be rank deficient by Corollary \ref{cor:rank_drop}. On the other hand, any such form has corank at most 1 by Corollary \ref{cor:corank1}, and hence its rank equals $n - 4$. Thus the rank is constant over all forms in ${\cal F}_{\bf u}' \setminus \{0\}$. Moreover, inequalities \eqref{ineq_relations} are satisfied strictly for such forms. Hence ${\cal F}_{\bf u}' \setminus \{0\}$ must be contained in the relative interior of ${\cal F}_{\bf u}'$, which implies that ${\cal F}_{\bf u}'$ is a ray generated by a single form $B'$. By Theorem \ref{thm:sdp_rep} $A' = \Lambda^{-1}(B')$ is then an exceptional extremal copositive matrix, and $F_{\bf u} \simeq \mathbb R_+^2$ is generated by $A'$ and $A_P$. Since $\rk B' = n-4$, the submatrices $A_{I_j}'$ also have rank $n-4$ by Corollary \ref{cor:rank_equal}, and hence corank 2.
\end{proof}

{\lemma \label{full_rank_odd} Let $n \geq 5$ be odd, and suppose $r_{\max} = n - 3$ and ${\cal F}_{\bf u} \not= {\cal P}_{\bf u}$. Then $F_{\bf u}$ does not contain non-zero positive semi-definite matrices.

If $F_{\bf u}$ is 1-dimensional, then it is generated by an extremal exceptional copositive matrix $A$ such that the submatrices $A_{I_j}$ have corank 1 for all $j = 1,\dots,n$.

If $\dim F_{\bf u} > 1$, then the monodromy $\mathfrak M$ of the system ${\bf S}_{\bf u}$ possesses the eigenvalue $-1$, and all boundary rays of $F_{\bf u}$ are generated by extremal exceptional copositive matrices. For any such boundary matrix $A$, its submatrices $A_{I_j}$ have corank 2 for all $j = 1,\dots,n$. }

\begin{proof}
As in the proof of the previous lemma, $\mathfrak M$ has all eigenvalues on the unit circle, with equal geometric and algebraic multiplicities. However, now $\dim {\cal L}_{\bf u}^* = n-3$ is even, and by \eqref{det_positive} the real eigenvalues $\pm1$, if they appear, have even multiplicity. By Corollaries \ref{cor:eigenvalue1} and \ref{cor:difference_ranks} the multiplicity of the eigenvalue 1 cannot exceed 1 and this eigenvalue does not appear. By Lemma \ref{eigenvalue1} we get ${\cal P}_{\bf u} = \{0\}$ and therefore by Theorem \ref{thm:sdp_rep} the face $F_{\bf u}$ does not contain non-zero positive semi-definite matrices.

Suppose now that $F_{\bf u}$ is not 1-dimensional. Then $\partial {\cal F}_{\bf u} \setminus {\cal P}_{\bf u} = \partial {\cal F}_{\bf u} \setminus \{0\}$ is not empty and by Corollary \ref{cor:rank_drop} consists of rank deficient forms. On the other hand, the rank can drop at most by 1 by virtue of Corollary \ref{cor:corank1}. Hence $\dim\ker B = 1$ for all $B \in \partial {\cal F}_{\bf u} \setminus \{0\}$. This kernel must be a real eigenspace of $\mathfrak M^*$, because $\mathfrak M^*$ preserves the form $B$ by shift-invariance. Therefore $\mathfrak M^*$ must have a real eigenvalue, which can only be equal to $-1$. Moreover, since the boundary subset $\partial {\cal F}_{\bf u} \setminus \{0\}$ consists of forms of corank 1, it must be smooth and every boundary ray is extremal. Indeed, if $\partial {\cal F}_{\bf u}$ would contain a face of dimension exceeding 1, then the boundary of this face would be different from $\{0\}$, but no further rank drop could occur there. Hence the pre-image $\Lambda^{-1}[\partial {\cal F}_{\bf u} \setminus \{0\}]$ consists of extremal exceptional copositive matrices. By Corollary \ref{cor:rank_equal} the submatrices $A_{I_j}$ of such a matrix $A$ have rank $n-4$, or corank 2, for all $j = 1,\dots,n$.

If, on the contrary, $\dim F_{\bf u} = 1$, then $F_{\bf u}$ is generated by an extremal exceptional copositive matrix $A$. By assumption the form $B = \Lambda(A)$ is positive definite and has rank $n-3$. By Corollary \ref{cor:rank_equal} the submatrices $A_{I_j}$ also have rank $n-3$ for all $j = 1,\dots,n$.
\end{proof}


Finally, we shall investigate the zero set of exceptional copositive matrices $A \in F_{\bf u}$.

{\lemma \label{zeros_in_Ij} Suppose $n \geq 5$, and let $A \in F_{\bf u}$ be an exceptional copositive matrix. If $v \in \mathbb R_+^n$ is a zero of $A$, then there exists an index $j \in \{1,\dots,n\}$ such that $\Supp v \subset I_j$. }

\begin{proof}
Since $A$ is exceptional, we cannot have $\Supp v = \{1,\dots,n\}$ by Lemma \ref{zeroPSD}. Therefore $\Supp v$ is a subset of $\{1,\dots,n\} \setminus \{k\}$ for some $k \in \{1,\dots,n\}$. Without loss of generality, assume that $\Supp v \subset \{1,\dots,n-1\}$. By Theorem \ref{first_cond} we have $(u^1)^TAu^2 > 0$, and hence by Lemma \ref{completion} we have $A_{(1,\dots,n-1)} = P + \delta E_{1,n-1}$ for some positive semi-definite matrix $P \in {\cal S}_+^{n-1}$ and some $\delta > 0$. Since $v$ is a zero of $A$, we get that the subvector $v_{(1,\dots,n-1)} \in \mathbb R_+^{n-1}$ is a zero of both $P$ and $E_{1,n-1}$. It follows that the first and the last element of this subvector cannot be simultaneously positive, which implies $\Supp v \subset I_1$ or $\Supp v \subset I_2$.
\end{proof}

{\theorem \label{thm:reg_deg} Let $A \in F_{\bf u}$ be an exceptional copositive matrix and set $B = \Lambda(A)$. Then either
\begin{itemize}
\setlength{\itemsep}{1pt}
\setlength{\parskip}{0pt}
\setlength{\parsep}{0pt}
\item[(i.a)] $A$ has minimal circulant zero support set;
\item[(i.b)] $B$ is positive definite;
\item[(i.c)] the corank of the submatrices $A_{I_j}$ equals 1, $j = 1,\dots,n$;
\item[(i.d)] the minimal zero support set of $A$ is $\{I_1,\dots,I_n\}$, with minimal zeros $u^1,\dots,u^n$;
\item[(i.e)] for even $n$ the matrix $A$ is the sum of an exceptional copositive matrix with non-minimal circulant zero support set and a rank 1 positive semi-definite matrix;
\item[(i.f)] if $n$ is odd and the monodromy operator $\mathfrak M$ has no eigenvalue equal to $-1$, then $A$ is extremal;
\end{itemize}
or
\begin{itemize}
\setlength{\itemsep}{1pt}
\setlength{\parskip}{0pt}
\setlength{\parsep}{0pt}
\item[(ii.a)] $A$ has non-minimal circulant zero support set;
\item[(ii.b)] the corank of $B$ equals 1;
\item[(ii.c)] the corank of the submatrices $A_{I_j}$ equals 2, $j = 1,\dots,n$;
\item[(ii.d)] the support of any minimal zero of $A$ is a strict subset of one of the index sets $I_1,\dots,I_n$, and every index set $I_j$ has exactly two subsets which are supports of minimal zeros of $A$;
\item[(ii.e)] every non-minimal zero of $A$ has support equal to $I_j$ for some $j = 1,\dots,n$ and is a sum of two minimal zeros;
\item[(ii.f)] $A$ is extremal.
\end{itemize} }

\begin{proof}
By Corollary \ref{cor:corank1} the form $B$ is either positive definite or has corank 1. By Corollary \ref{cor:rank_equal} the submatrices $A_{I_j}$ have corank 1 in the first case and corank 2 in the second case, for all $j = 1,\dots,n$. Hence by Lemma \ref{minimal_corank1} the zeros $u^j$ are minimal in the first case and not minimal in the second case, for all $j = 1,\dots,n$. By Lemma \ref{zeros_in_Ij} $A$ has minimal circulant zero support set in the first case, and there are no index sets other than $I_1,\dots,I_n$ which are supports of minimal zeros of $A$. Thus either (i.a)---(i.d) or (ii.a)---(ii.c) hold.

Assume the second case. By Lemma \ref{zeros_in_Ij} every support of a minimal zero of $A$ is a subset of one of the sets $I_1,\dots,I_n$. By Lemma \ref{minimal_corank1} this inclusion is strict. The set of zeros $v$ of $A$ satisfying $\Supp v \subset I_j$ is determined by the intersection of the two-dimensional kernel of $A_{I_j}$ with the nonnegative orthant. However, every two-dimensional convex cone is linearly isomorphic to $\mathbb R_+^2$ and is hence the convex hull of two extreme rays. Assertions (ii.d),(ii.e) then readily follow. The maximal rank achieved by forms in ${\cal F}_{\bf u}$ equals either $n-4$ or $n-3$ by Corollary \ref{cor:difference_ranks}. In the first case (ii.f) follows from Lemma \ref{exc_degenerate}, in the second case from Lemma \ref{full_rank_even} or \ref{full_rank_odd}, dependent on the parity of $n$.

Assume (i.a)---(i.d). Then (i.e) and (i.f) follow from Lemmas \ref{full_rank_even} and \ref{full_rank_odd}, respectively.
\end{proof}

The prototype of exceptional copositive matrices satisfying conditions (i.a)---(i.f) are the Hildebrand matrices \eqref{T_matrices}, while the prototype of those satisfying (ii.a)---(ii.f) is the Horn matrix \eqref{Horn_matrix}.

\section{Submanifolds of extremal exceptional copositive matrices}

In the previous two sections we considered the face $F_{\bf u} \subset {\cal C}^n$ for a fixed collection $\bf u$ of zeros. In Theorem \ref{thm:reg_deg} we have shown that there are two potential possibilities for an exceptional copositive matrix $A$ in such a face $F_{\bf u}$. Namely, either $A$ has minimal, or $A$ has non-minimal circulant zero support set, either case imposing its own set of conditions on $A$. In this section we show that in each of these cases, the matrix $A$ is embedded in a submanifold of ${\cal S}^n$ of codimension $n$ or $2n$, respectively, which consists of exceptional copositive matrices with similar properties. However, different matrices in this submanifold may belong to faces $F_{\bf u}$ corresponding to different collections $\bf u$.

{\theorem \label{thm:manifold_n} Let $n \geq 5$, and let $\hat A \in {\cal C}^n$ be an exceptional matrix with minimal circulant zero support set and with zeros $\hat u^1,\dots,\hat u^n \in \mathbb R_+^n$ such that $\Supp \hat u^j = I_j$. Then there exists a neighbourhood ${\cal U} \subset {\cal S}^n$ of $\hat A$ with the following properties:
\begin{itemize}
\setlength{\itemsep}{1pt}
\setlength{\parskip}{0pt}
\setlength{\parsep}{0pt}
\item[(i)] if $A \in {\cal U}$ and $\det A_{I_j} = 0$ for all $j = 1,\dots,n$, then $A$ is an exceptional copositive matrix with minimal circulant zero support set;
\item[(ii)] the set of matrices $A \in {\cal U}$ satisfying the conditions in (i) is an algebraic submanifold of codimension $n$ in ${\cal S}^n$.
\end{itemize} }

\begin{proof}
By Theorem \ref{thm:reg_deg} the submatrices $\hat A_{I_j}$ have rank $n-3$ for $j = 1,\dots,n$. Let $A \in {\cal S}^n$ be sufficiently close to $\hat A$ and such that $\det A_{I_j} = 0$ for all $j = 1,\dots,n$. For every $j = 1,\dots,n$ the submatrix $A_{I_j}$ has $n-3$ positive eigenvalues by continuity and one zero eigenvalue by assumption, and hence is positive semi-definite. Moreover, for every $j = 1,\dots,n$ the kernel of $A_{I_j}$ is close to the kernel of $\hat A_{I_j}$ and hence is generated by an element-wise positive vector. We then find vectors $u^1,\dots,u^n \in \mathbb R_+^n$, close to $\hat u^1,\dots,\hat u^n$, such that $\Supp u^j = I_j$ and the subvector $u^j_{I_j}$ is in the kernel of $A_{I_j}$. Then $(u^j)^TAu^j = 0$ for all $j = 1,\dots,n$, and $(u^j)^TAu^k$ is close to $(\hat u^j)^T\hat A\hat u^k$ for all $j,k = 1,\dots,n$. By Theorem \ref{first_cond} we have $(\hat u^n)^T\hat A\hat u^1 > 0$, $(\hat u^j)^T\hat A\hat u^{j+1} > 0$, $j = 1,\dots,n-1$, and hence also $(u^n)^TAu^1 > 0$, $(u^j)^TAu^{j+1} > 0$, $j = 1,\dots,n-1$. By Theorem \ref{first_cond} we then get that $A$ is a copositive exceptional matrix. By Theorem \ref{thm:reg_deg} the matrix $A$ has minimal circulant zero support set, which proves (i).

Consider the set ${\cal M} = \{ X \in {\cal S}^n \,|\, \det X_{I_j} = 0\ \forall\ j = 1,\dots,n \}$. It is defined by $n$ polynomial equations, and $\hat A \in {\cal M}$. By virtue of Lemma \ref{corank1_manifold} the gradient of the function $\det X_{I_j}$ at $X = \hat A$ is proportional to the rank 1 matrix $\hat u^j(\hat u^j)^T$. By Lemmas \ref{periodic_char}, \ref{eigenvalue1}, and Corollary \ref{cor:difference_ranks} the linear span of the zeros $\hat u^j$ has dimension at least $n-1$. Since no two of these zeros are proportional, the gradients of the functions $\det X_{I_j}$ are linearly independent at $X = \hat A$ by Lemma \ref{lem:gradients}. It follows that ${\cal M}$ is a smooth algebraic submanifold of codimension $n$ in a neighbourhood of $\hat A$. This proves (ii).
\end{proof}

{\theorem \label{thm:manifold_n_ext} The extremal exceptional copositive matrices with minimal circulant zero support set form an open subset of the manifold of all exceptional matrices with minimal circulant zero support set. }

\begin{proof}
Let $\hat A \in {\cal C}^n$ be as in the conditions of the previous theorem, but suppose in addition that $\hat A$ is extremal. Assume the notations in the proof of the previous theorem. The inequalities $(u^n)^TAu^1 > 0$, $(u^j)^TAu^{j+1} > 0$, $j = 1,\dots,n-1$ imply that the element $(Au^j)_i$ vanishes if and only if $i \in I_j$. Thus by Lemma \ref{extreme_criterion} the matrix $A$ is extremal if and only if the linear system of equations in the matrix $X \in {\cal S}^n$ given by
\begin{equation} \label{system_extreme}
(Xu^j)_i = 0\qquad \forall\ i \in I_j,\ j = 1,\dots,n
\end{equation}
has a 1-dimensional solution space. The coefficients of this system depend linearly on the entries of the zeros $u^j$. Moreover, since $\hat A$ is extremal, system \eqref{system_extreme} has a 1-dimensional solution space for $u^j = \hat u^j$. But the rank of a matrix is a lower semi-continuous function, hence the solution space of \eqref{system_extreme} has dimension at most 1 if the zeros $u^j$ are sufficiently close to $\hat u^j$. However, $X = A$ is always a solution, and therefore $A$ is an extremal copositive matrix.
\end{proof}

The simplest manifold of the type described in Theorem \ref{thm:manifold_n_ext} is the 10-dimensional union of the ${\cal G}_5$-orbits of the Hildebrand matrices \eqref{T_matrices}. Matrices \eqref{T_matrices} themselves depend on 5 parameters, while the action of ${\cal G}_5$ adds another 5 parameters.

{\theorem \label{thm:manifold_2n} Let $n \geq 5$, and let $\hat A \in {\cal C}^n$ be an exceptional matrix with non-minimal circulant zero support set and having zeros $\hat u^1,\dots,\hat u^n \in \mathbb R_+^n$ such that $\Supp \hat u^j = I_j$. Then there exists a neighbourhood ${\cal U} \subset {\cal S}^n$ of $\hat A$ with the following properties:
\begin{itemize}
\setlength{\itemsep}{1pt}
\setlength{\parskip}{0pt}
\setlength{\parsep}{0pt}
\item[(i)] if $A \in {\cal U}$ and $\rk A_{I_j} = n-4$ for all $j = 1,\dots,n$, then $A$ is an exceptional extremal copositive matrix with non-minimal circulant zero support set;
\item[(ii)] the set of matrices $A \in {\cal U}$ satisfying the conditions in (i) is an algebraic submanifold of codimension $2n$ in ${\cal S}^n$.
\end{itemize} }

\begin{proof}
By Theorem \ref{thm:reg_deg} the submatrices $\hat A_{I_j}$ have rank $n-4$ for $j = 1,\dots,n$. Let $A$ be sufficiently close to $\hat A$ and suppose that $\rk A_{I_j} = n-4$ for all $j = 1,\dots,n$. By continuity the $n-4$ largest eigenvalues of $A_{I_1}$ are then positive. However, the remaining two eigenvalues of $A_{I_1}$ are zero by assumption, and hence $A_{I_1} \succeq 0$. Now the kernel of $A_{I_1}$ is close to that of $\hat A_{I_1}$, hence $\ker A_{I_1}$ contains a positive vector close to $\hat u^1_{I_1}$. Then there exists a vector $u^1 \in \mathbb R_+^n$, close to $\hat u^1$, such that $\Supp u^1 = I_1$ and $(u^1)^TAu^1 = 0$. In a similar way we construct vectors $u^j \in \mathbb R_+^n$, close to $\hat u^j$, such that $\Supp u^j = I_j$ and $(u^j)^TAu^j = 0$, for all $j = 1,\dots,n$. By virtue of Theorem \ref{first_cond} we have $(\hat u^n)^T\hat A\hat u^1 > 0$ and $(\hat u^j)^T\hat A\hat u^{j+1} > 0$ for all $j = 1,\dots,n-1$. By continuity we then get $(u^n)^TAu^1 > 0$ and $(u^j)^TAu^{j+1} > 0$ for all $j = 1,\dots,n-1$. Again by Theorem \ref{first_cond} it then follows that $A$ is an exceptional copositive matrix. By Theorem \ref{thm:reg_deg} the matrix $A$ is also extremal and has non-minimal circulant zero support set, which proves (i).

The proof of (ii) is a bit more complicated. Set $\hat {\bf u} = \{\hat u^1,\dots,\hat u^n\}$. By Lemma \ref{lem:second_cond} we have $\hat A \in {\cal A}_{\hat{\bf u}}$. Let $\hat B = \Lambda(\hat A)$ be the positive semi-definite symmetric bilinear form corresponding to $\hat A$.

By Corollary \ref{cor:rank_equal} each of the principal submatrices $\hat A_{I_j'}$ has a 1-dimensional kernel. For every $j = 1,\dots,n$, we define a vector $\hat v^j \in \mathbb R^n$ such that $\hat v^j_i = 0$ for all $i \not\in I_j'$ and the subvector $\hat v^j_{I_j'}$ generates the kernel of $\hat A_{I_j'}$. We now claim the following:
\begin{itemize}
\setlength{\itemsep}{1pt}
\setlength{\parskip}{0pt}
\setlength{\parsep}{0pt}
\item[(a)] for every $j = 1,\dots,n-1$ we have $\ker \hat A_{I_j} = \spa\{ \hat v^j,\hat v^{j+1} \}$, and $\ker \hat A_{I_n} = \spa\{ \hat v^n,\hat v^1 \}$;
\item[(b)] the first and the last element of the subvector $\hat v^j_{I_j'}$ is non-zero for all $j = 1,\dots,n$;
\item[(c)] the vectors $\hat v^1,\dots,\hat v^n$ are linearly independent.
\end{itemize}

Since the positive semi-definite quadratic form $\hat A_{I_2}$ vanishes on the subvectors $\hat v^2_{I_2},\hat v^3_{I_2}$, we have that $\hat v^2_{I_2},\hat v^3_{I_2} \in \ker \hat A_{I_2}$. Now $\dim\ker\hat A_{I_2} = 2$, and hence either $\ker\hat A_{I_2} = \spa\{ \hat v^2_{I_2},\hat v^3_{I_2} \}$, or $\hat v^2,\hat v^3$ are linearly dependent.

Assume the latter for the sake of contradiction. Then the non-zero elements of $\hat v^2$ have indices in the intersection $I_2' \cap I_3' = \{3,\dots,n-2\}$. The relation $\hat A_{I_2}\hat v^2_{I_2} = 0$ together with \eqref{band_relation} yields $\sum_{s=3}^{n-2} \hat v^2_s\hat B({\bf e}_t,{\bf e}_s) = 0$ for all $t = 2,\dots,n-1$. The evaluation functionals ${\bf e}_2,\dots,{\bf e}_{n-1}$ span the whole space ${\cal L}_{\hat{\bf u}}^*$, and therefore we must have $\sum_{s=3}^{n-2} \hat v^2_s\hat B({\bf e}_t,{\bf e}_s) = 0$ for all $t \geq 1$. In particular, we get $\sum_{s=3}^{n-2} \hat v^2_s\hat B({\bf e}_t-{\bf e}_{t+n},{\bf e}_s) = 0$ for all $t = 1,\dots,n-4$. Since $\hat A$ is exceptional, by virtue of Lemma \ref{lem:triangular} the coefficient matrix of this linear homogeneous system on the elements of $\hat v^2$ is regular. Therefore $\hat v^2 = 0$, leading to a contradiction.

It follows that $\ker\hat A_{I_2} = \spa\{ \hat v^2_{I_2},\hat v^3_{I_2} \}$, and by repeating the argument after circular shifts of the indices we obtain (a). Since $\hat u^1 \in \ker\hat A_{I_1}$, we have $\hat u^1 = \alpha \hat v^1 + \beta \hat v^2$ for some coefficients $\alpha,\beta$. In particular, we have $\hat u^1_1 = \alpha\hat v^1_1 > 0$, implying $\hat v^1_1 \not = 0$. Similarly, $\hat u^1_{n-2} = \beta\hat v^2_{n-2} > 0$ implies $\hat v^2_{n-2} \not= 0$. Repeating the argument after circular shifts of the indices we obtain (b).

Finally, let $w \in \mathbb R^n$ be such that $\langle \hat v^j,w \rangle = 0$ for all $j = 1,\dots,n$. By (a) and by Lemma \ref{zeros_in_Ij} we have that every zero of $\hat A$ is a linear combination of the vectors $\hat v^j$. It follows that $w$ is orthogonal to all zeros of $\hat A$. By Lemma \ref{irred_rank1} there exists $\varepsilon > 0$ such that $\hat A - \varepsilon ww^T \in {\cal C}^n$. But by Theorem \ref{thm:reg_deg} $\hat A$ is extremal. Therefore we must have $w = 0$, which yields (c).

We now prove (ii). Consider the set
\[ {\cal M} = \left\{ X \in {\cal S}^n \,\left|\, \begin{array}{ll} \det X_{I'_j} = 0\quad \forall j = 1,\dots,n; \\ \det(X_{ik})_{i \in I'_j;k \in I'_{j+1}} = 0\quad \forall j = 1,\dots,n-1;\ \det(X_{ik})_{i \in I'_n;k \in I'_1} = 0 \end{array} \right. \right\}.
\]
It is defined by $2n$ polynomial equations, and $\hat A \in {\cal M}$. By Lemma \ref{corank2_manifold} in a neighbourhood of $\hat A$ the manifold ${\cal M}$ coincides with the set of matrices $A$ satisfying condition (i).

Now by virtue of Lemma \ref{corank1_manifold} the gradient of the function $\det X_{I'_j}$ at $X = \hat A$ is proportional to the rank 1 matrix $\hat v^j(\hat v^j)^T$, the gradient of $\det(X_{ik})_{i \in I'_j;k \in I'_{j+1}}$ is proportional to the rank 1 matrix $\hat v^j(\hat v^{j+1})^T$, and the gradient of $\det(X_{ik})_{i \in I'_n;k \in I'_1}$ is proportional to $\hat v^n(\hat v^1)^T$. By (c) the vectors $\hat v^j$ are linearly independent, hence these $2n$ gradients are also linearly independent. It follows that ${\cal M}$ is a smooth algebraic submanifold of codimension $2n$ in a neighbourhood of $\hat A$. This proves (ii).
\end{proof}

The simplest manifold of the type described in Theorem \ref{thm:manifold_2n} is the 5-dimensional ${\cal G}_5$-orbit of the Horn matrix \eqref{Horn_matrix}.

\section{Existence of non-trivial faces} \label{sec:existence}

So far we have always supposed that the feasible sets ${\cal F}_{\bf u}$ or ${\cal P}_{\bf u}$ of the LMIs in Theorem \ref{thm:sdp_rep} contain non-zero forms. In this section we shall explicitly construct non-zero faces ${\cal F}_{\bf u}$ and ${\cal P}_{\bf u}$ for arbitrary matrix sizes $n \geq 5$.

\subsection{Faces consisting of positive semi-definite matrices}

In this subsection we construct non-zero faces $F_{\bf u}$ of ${\cal C}^n$ which contain only positive semi-definite matrices, i.e., which satisfy $F_{\bf u} = P_{\bf u}$.

We shall need the following concept of a {\it slack matrix}, which has been introduced in \cite{Yannakakis91} for convex polytopes. Let $K \subset \mathbb R^m$ be a polyhedral convex cone, and let $K^* = \{ f \in \mathbb R_m \,|\, \langle f,x \rangle \geq 0\ \forall\ x \in K \}$ be its dual cone, where $\mathbb R_m$ is the space of linear functionals on $\mathbb R^m$. Then $K^*$ is also a convex polyhedral cone. Let $x_1,\dots,x_r$ be generators of the extreme rays of $K$, and $f_1,\dots,f_s$ generators of the extreme rays of $K^*$.

{\definition Assume the notations of the previous paragraph. The {\it slack matrix} of $K$ is the nonnegative $s \times r$ matrix $(\langle f_i,x_j \rangle)_{i = 1,\dots,s;j = 1,\dots,r}$. }

{\theorem \label{polygon_full_periodic} Assume $n \geq 5$, and let ${\bf u} = \{u^1,\dots,u^n\} \subset \mathbb R_+^n$ be such that $\Supp u^j = I_j$ for all $j = 1,\dots,n$. Let $U$ be the $n \times n$ matrix with columns $u^1,\dots,u^n$. Then the face $F_{\bf u}$ consists of positive semi-definite matrices up to rank $n-3$ inclusive if and only if $U$ is the slack matrix of a convex polyhedral cone $K \subset \mathbb R^3$ with $n$ extreme rays. }

\begin{proof}
By Theorem \ref{thm:sdp_rep} and Lemma \ref{PSD_full_rank} we have $F_{\bf u} = P_{\bf u} \simeq {\cal S}_+^{n-3}$ if and only if $\rk U = 3$.


Assume $\rk U = 3$. Choose a factorization $U = U_LU_R^T$, with $U_L,U_R$ being rank 3 matrices of size $n \times 3$. No two columns of $U$ and hence no two rows of $U_R$ are proportional, because $I_j \not\subset I_k$ for every $j \not= k$. Denote the convex conic hull of the rows of $U_R$ by $K$. Row $n$ of $U_L$ is orthogonal to rows 1 and 2 of $U_R$ and has a positive scalar product with the other rows of $U_R$. Therefore row $n$ of $U_L$ defines a supporting hyperplane to $K$, and the two-dimensional convex conic hull of row 1 and row 2 of $U_R$ is a subset of the boundary of $K$. By a circular shift of the indices and by repeating the argument we extend the construction of the boundary of $K$ until it closes in on itself. We obtain that $K$ is a polyhedral cone with $n$ extreme rays generated by the rows of $U_R$. On the other hand, the rows of $U_L$ define exactly the supporting hyperplanes to $K$ which intersect $K$ in its two-dimensional faces. Therefore the dual cone $K^*$ is given by the convex conic hull of the rows of $U_L$. The relation $U = U_LU_R^T$ then reveals that $U$ is the slack matrix of $K$.

On the other hand, a convex polyhedral cone $K \subset \mathbb R^3$ with $n \geq 3$ extreme rays as well as its dual $K^*$ have a linear span of dimension 3. Therefore the slack matrix of such a cone has rank 3. This completes the proof.
\end{proof}

Theorem \ref{polygon_full_periodic} provides a way to construct all collections ${\bf u} \subset \mathbb R_+^n$ of vectors $u^1,\dots,u^n$ satisfying $\Supp u^j = I_j$, $j = 1,\dots,n$, such that the face $F_{\bf u}$ of ${\cal C}^n$ consists of positive semi-definite matrices only and is linearly isomorphic to ${\cal S}_+^{n-3}$.

Let now ${\bf u} = \{u^1,\dots,u^n\} \subset \mathbb R_+^n$ be an arbitrary collection such that $\Supp u^j = I_j$ for all $j = 1,\dots,n$. Let again $U$ be the $n \times n$ matrix with columns $u^1,\dots,u^n$. By Lemmas \ref{periodic_char} and \ref{eigenvalue1} we have ${\cal P}_{\bf u} \simeq {\cal S}^k_+$, where $k$ is the corank of $U$. We have shown that there exist collections ${\bf u}$ such that $U$ has corank $n-3$. By perturbing some of the zeros $u^j$ in such a collection, the corank of $U$ can be decreased and may assume an arbitrary value between 0 and $n-3$. In this way we obtain faces $F_{\bf u}$ of ${\cal C}^n$ for which the subset $P_{\bf u}$ of positive semi-definite matrices is isomorphic to ${\cal S}_+^k$ with arbitrary rank $k = 0,\dots,n-3$. For $k \geq 2$ we have by Corollary \ref{cor:difference_ranks} that $F_{\bf u} = P_{\bf u}$.

\subsection{Circulant matrices} \label{subs:circulant}

In this section we consider faces $F_{\bf u}$ defined by special collections $\bf u$. Let $u \in \mathbb R_{++}^{n-2}$ be {\it palindromic}, i.e., with positive entries and invariant with respect to inversion of the order of its entries. Define ${\bf u} = \{ u^1,\dots,u^n \} \subset \mathbb R_+^n$ such that $\Supp u^j = I_j$ and $u^j_{I_j} = u$ for all $j = 1,\dots,n$. By construction, the linear dynamical system ${\bf S}_{\bf u}$ defined by $\bf u$ has constant coefficients, namely the entries of $u$. Set $p(x) = \sum_{k=0}^{n-3} u_{k+1}x^k$.

We provide necessary and sufficient conditions on $\bf u$ such that the corresponding face $F_{\bf u} \subset {\cal C}^n$ contains exceptional copositive matrices, and construct explicit collections $\bf u$ which satisfy these conditions. We show that the copositive matrices in these faces must be circulant, i.e., invariant with respect to simultaneous circular shifts of its row and column indices.


{\lemma \label{lem:circ_ex} Suppose the collection $\bf u$ is as in the first paragraph of this section. If $F_{\bf u} \not= P_{\bf u}$, then $F_{\bf u}$ contains exceptional circulant matrices. }

\begin{proof}
Let $A^1 \in F_{\bf u}$ be exceptional. By repeated simultaneous circular shifts of the row and column indices of $A^1$ by one entry we obtain copositive matrices $A^2,\dots,A^n$ which are also elements of $F_{\bf u}$. Then $A = \frac1n\sum_{j=1}^n A^j \in F_{\bf u}$ is a copositive circulant matrix. By Theorem \ref{first_cond} it is also exceptional.
\end{proof}

{\lemma \label{lem:Btoep} Suppose the collection $\bf u$ is as above, let $n \geq 5$, and let $A \in F_{\bf u}$ be a circulant matrix. Then the symmetric bilinear form $B = \Lambda(A)$ satisfies $B({\bf e}_t,{\bf e}_s) = B({\bf e}_{t+l},{\bf e}_{s+l})$ for all $t,s,l \geq 1$. }

\begin{proof}
Define matrices $B_{r,l} = (B({\bf e}_{r+t},{\bf e}_{r+s}))_{t,s = 0,\dots,l}$ of size $(l+1) \times (l+1)$, $r \geq 1$, $l \geq d$. By Lemma \ref{lem:imageL} for $l = d$ the matrices $B_{r,l}$ are equal to a single matrix $B_l$ for all $r \geq 1$, and this positive semi-definite matrix $B_l$ is Toeplitz, rank-deficient, and has an element-wise positive kernel vector. We now show by induction that this holds also for all $l > d$.

Indeed, the entries of the positive semi-definite matrices $B_{r,l+1}$ are all determined by the matrix $B_l$, except the upper right (and lower left) corner element. Since $B_l$ has a kernel vector with positive elements, Lemma \ref{completion} is applicable and these corner elements are unique and hence all equal for all $r \geq 1$. It follows that the matrices $B_{r,l+1}$ are all equal to a single positive semi-definite matrix $B_{l+1}$. By construction this matrix is also Toeplitz and rank-deficient, and by Lemma \ref{completion} it has a kernel vector with positive elements.

The claim of the lemma now easily follows.
\end{proof}

{\corollary \label{cor:trigB} Assume the conditions of the previous lemma and set $r = \lfloor \frac{\rk B}{2} \rfloor$. Then there exist distinct angles $\varphi_0 = \pi$, $\varphi_1,\dots,\varphi_r \in (0,\pi)$ and positive numbers $\lambda_0,\dots,\lambda_r$ such that
\begin{equation} \label{claimed_formula}
B({\bf e}_t,{\bf e}_s) = \left\{ \begin{array}{rcl} \sum_{j=0}^r \lambda_j \cos(|t-s|\varphi_j),&\quad& \rk B\ \mbox{odd}; \\
\sum_{j=1}^r \lambda_j \cos(|t-s|\varphi_j),&\quad& \rk B\ \mbox{even} \end{array} \right.
\end{equation}
for all $t,s \geq 1$. Moreover, $e^{\pm i\varphi_j}$ are roots of $p(x)$ for all $j = 1,\dots,r$. In addition, if $\rk B$ is odd, then $-1$ is also a root of $p(x)$. }

\begin{proof}
By Theorem \ref{thm:sdp_rep} we have $B \in {\cal F}_{\bf u}$. Hence the Toeplitz matrix $T = (B({\bf e}_i,{\bf e}_j))_{i,j = 1,\dots,n-2}$ is positive semi-definite, rank-deficient and of the same rank as $B$, and has the element-wise positive kernel vector $u$. Application of Lemma \ref{lem:Toeplitz_rk1sum} to $T$ then yields \eqref{claimed_formula} for $t,s = 1,\dots,N$ with $N = n-2$.

In order to show \eqref{claimed_formula} for $t,s = 1,\dots,N$ with arbitrary $N$ we shall proceed by induction over $N$. Suppose that \eqref{claimed_formula} is valid for $t,s = 1,\dots,N-1$, and that the matrix $T_{N-1} = (B({\bf e}_i,{\bf e}_j))_{i,j = 1,\dots,N-1}$ has a kernel vector $w_{N-1}$ with positive entries. The matrix $T_N = (B({\bf e}_i,{\bf e}_j)i)_{i,j = 1,\dots,N}$ is positive semi-definite by the positivity of $B$, and Toeplitz by Lemma \ref{lem:Btoep}. The second property implies that \eqref{claimed_formula} has to be proven only for the corner elements $(t,s) = (1,N)$ and $(t,s) = (N,1)$. Let us now consider all elements of $T_N$ except these corner elements as fixed (and given by \eqref{claimed_formula}) and the corner elements as variable. Positive semi-definiteness of $T_N$ implies that the value $B({\bf e}_1,{\bf e}_N)$ solves the corresponding positive semi-definite completion problem. By Lemma \ref{completion} the solution of this completion problem is unique. However, it is precisely the value given by \eqref{claimed_formula} which solves the completion problem. Therefore $B({\bf e}_1,{\bf e}_N)$ must also be given by \eqref{claimed_formula}.

Finally, by Lemma \ref{completion} $T_N$ has a kernel vector with positive entries, completing the induction step.
\end{proof}

Inserting \eqref{claimed_formula} into \eqref{linear_relation} for $t,s \geq 1$ such that $k = |t-s| = 3,\dots,\lceil \frac{n}{2} \rceil-1$ yields
\begin{equation} \label{linear_circulant}
\begin{array}{rclcl}
\sum_{j=0}^r \lambda_j (\cos((n-k)\varphi_j) - \cos(k\varphi_j)) &=& 0,&\quad& \rk B\ \mbox{odd}, \\
\sum_{j=1}^r \lambda_j (\cos((n-k)\varphi_j) - \cos(k\varphi_j)) &=& 0,&\quad& \rk B\ \mbox{even},
\end{array}\qquad 3 \leq k < \frac{n}{2}.
\end{equation}
Note that the range of $t,s$ is actually larger in \eqref{linear_relation}, but the relations are trivial for $k = \frac{n}{2}$ and yield \eqref{linear_circulant} again for $\frac{n}{2} < k \leq n-3$. Likewise, inserting \eqref{claimed_formula} into \eqref{ineq_relations} yields
\begin{equation} \label{ineq_circulant}
\begin{array}{rclcl}
\sum_{j=0}^r \lambda_j (\cos((n-2)\varphi_j) - \cos(2\varphi_j)) &\leq& 0,&\quad& \rk B\ \mbox{odd}, \\
\sum_{j=1}^r \lambda_j (\cos((n-2)\varphi_j) - \cos(2\varphi_j)) &\leq& 0,&\quad& \rk B\ \mbox{even}.
\end{array}
\end{equation}
Note that \eqref{linear_circulant} is a linear homogeneous system of equations in the weights $\lambda_j$.

{\lemma \label{lem:circulant_even} Let $n > 5$ be even, let $\bf u$ be as above, and let $A \in F_{\bf u}$ be an exceptional copositive circulant matrix. Then there exist $m = \frac{n}{2} - 2$ distinct angles $\zeta_1,\dots,\zeta_m \in (0,\pi)$, arranged in increasing order, with the following properties:

(a) the fractional part of $\frac{n\zeta_j}{4\pi}$ is in $(0,\frac12)$ for odd $j$ and in $(\frac12,1)$ for even $j$;

(b) the polynomial $p(x)$ is a positive multiple of $(x+1)\cdot\prod_{j=1}^m(x^2-2x\cos\zeta_j+1)$;

(c) there exist $c > 0$, $\lambda \geq 0$ such that for all $k = 1,\dots,\frac{n}{2}+1$ we have
\begin{equation} \label{formula_even}
A_{1k} = (-1)^{k-1}\lambda + c \cdot \sum_{j=1}^m \frac{\cos(k-1)\zeta_j}{\sin\zeta_j\sin\frac{n\zeta_j}{2}\prod_{l \not= j}(\cos\zeta_j - \cos\zeta_l)}.
\end{equation}
If $\lambda = 0$, then $A$ is extremal with non-minimal circulant zero support set. If $\lambda > 0$, then $A$ has minimal circulant zero support set and is not extremal. }

\begin{proof}
Let $B = \Lambda(A)$. By Corollary \ref{cor:corank1} we have either $\rk B = n - 4$ or $\rk B = n - 3$. Applying Corollary \ref{cor:trigB}, we have $r = \lfloor \frac{\rk B}{2} \rfloor = \frac{n}{2} - 2 = m$. Define $\zeta_1,\dots,\zeta_m$ to be the angles $\varphi_1,\dots,\varphi_r$, arranged in increasing order. Then $e^{\pm i\zeta_j}$ are roots of $p(x)$ by Corollary \ref{cor:trigB}, and since $(x - e^{i\zeta_j})(x - e^{-i\zeta_j}) = (x^2 - 2x\cos\zeta_j + 1)$, $p(x)$ is of the form $\alpha\cdot(x-\beta)\cdot\prod_{j=1}^m(x^2-2x\cos\zeta_j+1)$ for some $\alpha > 0$ and real $\beta$. Since the coefficient vector $u$ of $p(x)$ is palindromic, we must have $\beta = -1$, which proves (b).

By \eqref{band_relation} we have $A_{1k} = B({\bf e}_1,{\bf e}_k)$, $k = 1,\dots,\frac{n}{2}+1$, which in turn is given by Corollary \ref{cor:trigB}. Here the weights $\lambda_j$ satisfy relations \eqref{linear_circulant},\eqref{ineq_circulant}, the inequality being strict by Lemma \ref{lem:second_cond}. We have $\varphi_0 = \pi$ and hence $\cos((n-k)\varphi_0) = \cos(k\varphi_0)$ for all integers $k$. Therefore \eqref{linear_circulant},\eqref{ineq_circulant} do not impose any conditions on the coefficient $\lambda_0$, and these $m$ relations can be considered as conditions on $\lambda_1,\dots,\lambda_m$ only. By Corollary \ref{lambda_eq} there are no multiples of $\frac{2\pi}{n}$ among the angles $\zeta_1,\dots,\zeta_m$, and the coefficient vector $(\lambda_1,\dots,\lambda_m)$ is proportional to solution \eqref{solution_lambda} in this corollary, with a positive proportionality constant. This yields (c) with $c > 0$, with $\lambda = 0$ if $\rk B = n-4$, and with $\lambda = \lambda_0 > 0$ if $\rk B = n-3$, by virtue of Corollary \ref{cor:trigB}.

The last assertions of the lemma now follow from Theorem \ref{thm:reg_deg}.

Finally, the weights $\lambda_1,\dots,\lambda_m$ are positive by Corollary \ref{cor:trigB}, and Lemma \ref{coef_positive} implies (a).
\end{proof}

{\lemma \label{lem:circulant_odd} Let $n \geq 5$ be odd, let $\bf u$ be as above, and let $A \in F_{\bf u}$ be an exceptional copositive circulant matrix. Then there exist $m = \frac{n-3}{2}$ distinct angles $\zeta_1,\dots,\zeta_m \in (0,\pi]$, arranged in increasing order, with the following properties:

(a) the fractional part of $\frac{n\zeta_j}{4\pi}$ is in $(0,\frac12)$ for odd $j$ and in $(\frac12,1)$ for even $j$;

(b) the polynomial $p(x)$ is a positive multiple of $\prod_{j=1}^m(x^2-2x\cos\zeta_j+1)$;

(c) there exists $c > 0$ such that for all $k = 1,\dots,\frac{n+1}{2}$ we have
\begin{equation} \label{formula_odd}
A_{1k} = c \cdot \sum_{j=1}^m \frac{\cos(k-1)\zeta_j}{\sin\frac{\zeta_j}{2}\sin\frac{n\zeta_j}{2}\prod_{l \not= j}(\cos\zeta_j - \cos\zeta_l)}.
\end{equation}
The matrix $A$ has non-minimal circulant zero support set if $\zeta_m = \pi$ and minimal circulant zero support set if $\zeta_m < \pi$. In both cases $A$ is extremal. }

\begin{proof}
Let $B = \Lambda(A)$. By Corollary \ref{cor:corank1} we have either $\rk B = n-4$ or $\rk B = n-3$. Applying Corollary \ref{cor:trigB}, we get either $r = \lfloor \frac{\rk B}{2} \rfloor = \frac{n-5}{2} = m - 1$ or $r = \lfloor \frac{\rk B}{2} \rfloor = \frac{n-3}{2} = m$, respectively. Define $\zeta_1,\dots,\zeta_m$ to be the angles $\varphi_0,\dots,\varphi_r$ in the first case and $\varphi_1,\dots,\varphi_r$ in the second case, arranged in increasing order. Then $\zeta_m = \pi$ if $\rk B = n-4$ and $\zeta_m < \pi$ if $\rk B = n-3$. By Corollary \ref{cor:trigB} the numbers $e^{\pm i\zeta_j}$ are roots of $p(x)$, and $p(x)$ is of the form $\alpha\cdot(x+1)\cdot(x-\beta)\cdot\prod_{j=1}^{m-1}(x^2-2x\cos\zeta_j+1)$ if $\zeta_m = \pi$ and $p(x) = \alpha\cdot\prod_{j=1}^m(x^2-2x\cos\zeta_j+1)$ if $\zeta_m < \pi$, for some $\alpha > 0$ and real $\beta$. Since the coefficient vector $u$ of $p(x)$ is palindromic, we must have $\beta = -1$, which proves (b).

By Theorem \ref{thm:reg_deg} the matrix $A$ is extremal with non-minimal circulant zero support set if $\zeta_m = \pi$, and has minimal circulant zero support set if $\zeta_m < \pi$.

By \eqref{band_relation} we have $A_{1k} = B({\bf e}_1,{\bf e}_k)$, $k = 1,\dots,\frac{n+1}{2}$, which in turn is given by Corollary \ref{cor:trigB}. Here the $m$ coefficients $\lambda_j$ satisfy the $m$ relations \eqref{linear_circulant},\eqref{ineq_circulant}, the inequality being strict by Lemma \ref{lem:second_cond}. By Corollary \ref{lambda_eq} there are no multiples of $\frac{2\pi}{n}$ among the angles $\zeta_1,\dots,\zeta_m$, and the coefficient vector $(\lambda_1,\dots,\lambda_m)$ is proportional to solution \eqref{solution_lambda} in this corollary with a positive proportionality constant. This yields (c) with $c > 0$.

Assertion (a) is obtained in the same way as in the proof of Lemma \ref{lem:circulant_even}.

It remains to show the extremality of $A$ for $\zeta_m < \pi$. By Lemma \ref{lem:circ_extremal} in Appendix \ref{sec:appB} the dimension of the space ${\cal A}_{\bf u}$ cannot exceed 1, otherwise we obtain a contradiction with the nonnegativity of $u$. It follows that $\dim F_{\bf u} = 1$ and $A$ is extremal in ${\cal C}^n$.
\end{proof}

Note that the elements $A_{1k}$, $k = 1,\dots,\lceil \frac{n+1}{2} \rceil$, determine the matrix $A$ completely by its circulant property. We obtain the following characterization of collections $\bf u$ defining faces $F_{\bf u}$ which contain exceptional matrices.

{\theorem \label{thm:circ_char_even} Let $n > 5$ be even, $m = \frac{n}{2} - 2$, and let $\bf u$ and $p(x)$ be as in the first paragraph of this section. Then $F_{\bf u} \not= P_{\bf u}$ if and only if there exist distinct angles $\zeta_1,\dots,\zeta_m \in (0,\pi)$, arranged in increasing order, such that conditions (a),(b) of Lemma \ref{lem:circulant_even} hold. In this case the face $F_{\bf u}$ is linearly isomorphic to $\mathbb R_+^2$ and consists of the circulant matrices $A$ with entries $A_{1k}$, $k = 1,\dots,\frac{n}{2}+1$, given by \eqref{formula_even} with $c,\lambda \geq 0$. The subset $P_{\bf u} \subset F_{\bf u}$ is given by those $A$ with $c = 0$. }

\begin{proof}
Suppose $F_{\bf u} \not= P_{\bf u}$. Then by Lemma \ref{lem:circ_ex} there exists a circulant matrix $A \in F_{\bf u} \setminus P_{\bf u}$. Hence Lemma \ref{lem:circulant_even} applies and conditions (a),(b) hold.

Let now $\zeta_1,\dots,\zeta_m \in (0,\pi)$ be an increasing sequence of angles satisfying conditions (a),(b) of Lemma \ref{lem:circulant_even}. Define weights $\lambda_1,\dots,\lambda_m$ by \eqref{solution_lambda}. By Lemma \ref{coef_positive} these weights are positive, and by Corollary \ref{lambda_eq} they satisfy system \eqref{linear_inhom}. Let $A \in {\cal S}^n$ be the circulant matrix satisfying $A_{1k} = \sum_{j=1}^m \lambda_j\cos(k-1)\zeta_j$ for $k = 1,\dots,\frac{n}{2}+1$. By \eqref{linear_inhom} the last relation actually holds for $k = 1,\dots,n-2$, and $A_{1,n-1} > \sum_{j=1}^m \lambda_j\cos(n-2)\zeta_j$. By construction the submatrices $A_{I_j}$ are positive semi-definite of rank $2m = n-4$ and by condition (b) they possess the kernel vector $u = u^j_{I_j}$. Moreover, the preceding inequality implies that $(u^n)^TAu^1 > 0$ and $(u^j)^TAu^{j+1} > 0$ for $j = 1,\dots,n-1$. By Theorem \ref{first_cond} we then have $A \in F_{\bf u} \setminus P_{\bf u}$, proving the equivalence claimed in the theorem.

Let $P \in {\cal S}_+^n$ be the positive semi-definite circulant rank 1 matrix given element-wise by $P_{kl} = (-1)^{k-l}$. Since the subvectors $u^j_{I_j}$ are palindromic and have an even number of entries, we get that $(u^j)^TPu^j = 0$ for all $j = 1,\dots,n$. Hence $P \in P_{\bf u}$ and $r_{\it{PSD}} \geq 1$, where $r_{\it{PSD}}$ is the maximal rank achieved by matrices in $P_{\bf u}$. But then $r_{\it{PSD}} = 1$ and $r_{\max} = n-3$ by Corollary \ref{cor:difference_ranks}, and the last assertion of the theorem follows.

By Lemma \ref{full_rank_even} we have $F_{\bf u} \simeq \mathbb R_+^2$. However, the matrices $A$ and $P$ constructed above are linearly independent elements of $F_{\bf u}$, and therefore $F_{\bf u} \subset \spa\{A,P\}$ consists of circulant matrices only. The remaining assertions now follow from Lemma \ref{lem:circulant_even}.
\end{proof}

{\theorem \label{thm:circ_char_odd} Let $n \geq 5$ be odd, $m = \frac{n-3}{2}$, and let $\bf u$ and $p(x)$ be as in the first paragraph of this section. Then $F_{\bf u} \not= P_{\bf u}$ if and only if there exist distinct angles $\zeta_1,\dots,\zeta_m \in (0,\pi]$, arranged in increasing order, such that conditions (a),(b) of Lemma \ref{lem:circulant_odd} hold. In this case the face $F_{\bf u}$ is an extreme ray of ${\cal C}^n$ and consists of the circulant matrices $A$ with entries $A_{1k}$, $k = 1,\dots,\frac{n+1}{2}$, given by \eqref{formula_odd} with $c \geq 0$. }

\begin{proof}
The proof of the theorem is similar to the proof of Theorem \ref{thm:circ_char_even}, with obvious modifications.
%
%
\end{proof}

The question which collections $\bf u$, of the type described at the beginning of this subsection, yield faces $F_{\bf u}$ containing exceptional copositive matrices hence reduces to the characterization of real polynomials of the form given in (b) of Lemmas \ref{lem:circulant_even} or \ref{lem:circulant_odd}, with positive coefficients and satisfying condition (a) of these lemmas. This is seemingly a difficult question, and only limited results are known. However, the existence of faces $F_{\bf u}$ containing exceptional copositive matrices is guaranteed for every $n \geq 5$ by the following result on polynomials with equally spaced roots on the unit circle.

{\lemma \cite[Theorem 2]{EvansGreene91} \label{lem:Greene} Let $m \geq 1$ be an integer, and let $\alpha > 0$, $\theta \geq 0$ be such that $\frac{\pi}{2} \leq \theta + \frac{(m-1)\alpha}{2} \leq \pi$ and $0 < \alpha < \frac{\pi}{m}$. Then the polynomial $q(x) = \prod_{j=1}^m (x^2 - 2x\cos(\theta+(j-1)\alpha) + 1)$ has positive coefficients. }

\medskip

In order to construct explicit examples we need also the following result.

{\lemma \label{explicit_elements} Let $\tilde n \geq 5$ be an integer, and let $A$ be the $(\tilde n-2) \times (\tilde n-2)$ real symmetric Toeplitz matrix with first row $(2(1+2\cos\frac{\pi}{\tilde n}\cos\frac{3\pi}{\tilde n}), -2(\cos\frac{\pi}{\tilde n}+\cos\frac{3\pi}{\tilde n}), 1, 0, \dots, 0)$. Let further $m = \lceil \frac{\tilde n-4}{2} \rceil$, and set $\zeta_j = \frac{(2j+3)\pi}{\tilde n}$, $j = 1,\dots,m$. Then the following holds:
\begin{itemize}
\item if $\tilde n$ is even, then there exists $c > 0$ such that \eqref{formula_even} holds with $\lambda = 0$ and $n = \tilde n$ for all $k = 1,\dots,\frac{n}{2}+1$;
\item if $\tilde n$ is odd, then there exists $c > 0$ such that \eqref{formula_odd} holds with $n = \tilde n$ for all $k = 1,\dots,\frac{n+1}{2}$;
\item if $\tilde n$ is even, then there exists $c > 0$ such that \eqref{formula_odd} holds with $n = \tilde n - 1$ for all $k = 1,\dots,\frac{n+1}{2}$.
\end{itemize} } 

\begin{proof}
First note that $e^{\pm i\zeta_j}$ are roots of the polynomial $x^n + 1$ for all $j = 1,\dots,m$, along with $e^{\pm i\pi/n}$ and $e^{\pm 3i\pi/n}$ which are the 4 remaining roots. Hence $p(x) = \frac{x^n+1}{(x^2-2x\cos\frac{\pi}{n}+1)(x^2-2x\cos\frac{3\pi}{n}+1)}$ is a polynomial of degree $n-4$ whose roots are exactly $e^{\pm i\zeta_j}$, $j = 1,\dots,m$. Here for odd $n$ the root $e^{\pm i\zeta_m} = -1$ is counted only once. Let $w \in \mathbb R^{n-3}$ be the coefficient vector of $p(x)$.

Then the two vectors $u = (w^T, 0)^T$, $v = (0, w^T)^T$ are in the kernel of $A$. Indeed, $A$ is a Toeplitz band matrix with the elements of the band given by the coefficients of the denominator polynomial $d(x) = (x^2-2x\cos\frac{\pi}{n}+1)(x^2-2x\cos\frac{3\pi}{n}+1)$ in the above definition of $p(x)$. Multiplying a vector by $A$ hence amounts to a convolution of this vector with the coefficient vector of $d(x)$ and discarding the first two and the last two elements of the resulting vector. Therefore the elements of the vectors $Au$ and $Av$ are among those coefficients of the polynomial $x^n+1$ which do not correspond to the highest or the lowest power of $x$. It follows that $Au,Av$ vanish.

Let $A'$ be the upper left $(\tilde n-3) \times (\tilde n-3)$ principal submatrix of $A$. From the above it follows that $A'w = 0$. Then by Lemma \ref{lem:Toeplitz_linear_comb} the matrix $A'$ is a linear combination of the real symmetric Toeplitz matrices $T'_j$, $j = 1,\dots,m$, defined such that the first row of $T'_j$ equals $(1,\cos\zeta_j,\dots,\cos(n-4)\zeta_j)$. In other words, there exist real numbers $\mu_1,\dots,\mu_m$ such that $A' = \sum_{j=1}^m \mu_j T'_j$.

Let $T_j$, $j = 1,\dots,m$, be the real symmetric Toeplitz matrix with first row $(1,\cos\zeta_j,\dots,\cos(n-3)\zeta_j)$. Then the upper left and the lower right $(n-3) \times (n-3)$ principal submatrices of $T_j$ equal $T'_j$. Therefore all elements of the difference $D = A - \sum_{j=1}^m \mu_j T_j$ except possibly the upper right and the lower left corner element vanish. On the other hand, we have $Du = Dv = 0$, because $Au = Av = 0$ and $T_ju = T_jv = 0$ for all $j = 1,\dots,m$. But the first element of $u$ and the last element of $v$ equal 1, and thus the corner elements of $D$ vanish too.

Therefore $A = \sum_{j=1}^m \mu_j T_j$ and $A_{1k} = \sum_{j=1}^m \mu_j\cos(k-1)\zeta_j$ for all $k = 1,\dots,n-2$. On the other hand, we have $A_{1k} = 0$ for all $4 \leq k \leq n-2$ and hence $\sum_{j=1}^m \mu_j (\cos(n-k)\zeta_j - \cos k\zeta_j) = 0$ for all $3 \leq k \leq m + 1$. Set $c = -2^{-m}\sum_{j=1}^m \mu_j(\cos(n-2)\zeta_j - \cos2\zeta_j)$. Then the coefficients $\mu_j$ satisfy linear system \eqref{linear_inhom}, except that the right-hand side is multiplied by $2^m \cdot c$. Among the angles $\zeta_j$ there are no multiples of $\frac{2\pi}{n}$. By Corollary \ref{lambda_eq} we then have $\mu_j = 2^mc\lambda_j$, where $\lambda_j$ is given by \eqref{solution_lambda}.

This yields \eqref{formula_even} with $\lambda = 0$ for all $k = 1,\dots,\frac{n}{2}+1$ if $n$ is even, and \eqref{formula_odd} for all $k = 1,\dots,\frac{n+1}{2}$ if $n$ is odd. Plugging in $k = 1$, we obtain
\[ c = \frac{2(1+2\cos\frac{\pi}{n}\cos\frac{3\pi}{n})}{\sum_{j=1}^m \left(\sin\zeta_j\sin\frac{(2j+3)\pi}{2}\prod_{l \not= j}(\cos\zeta_j - \cos\zeta_l)\right)^{-1}}
\]
for even $n$ and
\[ c = \frac{2(1+2\cos\frac{\pi}{n}\cos\frac{3\pi}{n})}{\sum_{j=1}^m \left(\sin\frac{\zeta_j}{2}\sin\frac{(2j+3)\pi}{2}\prod_{l \not= j}(\cos\zeta_j - \cos\zeta_l)\right)^{-1}}
\]
for odd $n$. For all $j = 1,\dots,m$ we have $\sin\zeta_j > 0$ for even $n$ and $\sin\frac{\zeta_j}{2} > 0$ for odd $n$, $\sin\frac{(2j+3)\pi}{2} = (-1)^{j-1}$, $\prod_{l \not= j}(\cos\zeta_j - \cos\zeta_l) = (-1)^{j-1}$, and $1+2\cos\frac{\pi}{n}\cos\frac{3\pi}{n} > 0$. Hence $c > 0$, which completes the proof.
\end{proof}

These results allow to construct the following explicit examples.

{\sl Degenerate extremal matrices.} Let $n \geq 5$, $m = \lceil \frac{n}{2} \rceil - 2$, and $p(x) = \frac{(x^n+1)(x+1)}{(x^2-2x\cos\frac{\pi}{n}+1)(x^2-2x\cos\frac{3\pi}{n}+1)}$. Then $p(x)$ is a palindromic polynomial of degree $n-3$. Set also $q(x) = p(x)$ for odd $n$ and $q(x) = \frac{p(x)}{x+1}$ for even $n$. Then $q(x)$ is of degree $2m$ and has positive coefficients by virtue of Lemma \ref{lem:Greene} with $\alpha = \frac{2\pi}{n}$, $\theta = \frac{5\pi}{n}$. It follows that also $p(x)$ has positive coefficients. Let $u \in \mathbb R_+^{n-2}$ be the vector of its coefficients, and let $\bf u$ be the collection of nonnegative vectors constructed from $u$ as in the first paragraph of this section. Then the angles $\zeta_j = \frac{(2j+3)\pi}{n}$, $j = 1,\dots,m$, satisfy conditions (a),(b) of Lemmas \ref{lem:circulant_even} and \ref{lem:circulant_odd}, for even and odd $n$, respectively. By Theorems \ref{thm:circ_char_even} and \ref{thm:circ_char_odd} we obtain that $F_{\bf u} \simeq \mathbb R_+^2$ for even $n$ and $F_{\bf u} \simeq \mathbb R_+$ for odd $n$, their elements being circulant matrices given by \eqref{formula_even} and \eqref{formula_odd}, respectively. One extreme ray of $F_{\bf u}$ is then generated by an extremal copositive circulant matrix $A$ with non-minimal circulant zero support set. For even $n$ the other extreme ray is generated by a circulant positive semi-definite rank 1 matrix $P$. Their elements are given by $P_{ij} = (-1)^{i-j}$ and
\begin{equation} \label{example_deg}
A_{ij} = \left\{ \begin{array}{rcl} 2(1+2\cos\frac{\pi}{n}\cos\frac{3\pi}{n}),&\quad&i = j, \\
-2(\cos\frac{\pi}{n}+\cos\frac{3\pi}{n}),&&|i-j| \in \{1,n-1\}, \\
1,&&|i-j| \in \{2,n-2\}, \\
0,&&|i-j| \in \{3,\dots,n-3\}, \end{array} \right.
\end{equation}
$i,j = 1,\dots,n$ by the first two assertions of Lemma \ref{explicit_elements}.

\medskip

{\sl Regular extremal matrices.} Let $n \geq 5$ be odd, and set $m = \frac{n-3}{2}$, $p(x) = \frac{x^{n+1}+1}{(x^2-2x\cos\frac{\pi}{n+1}+1)(x^2-2x\cos\frac{3\pi}{n+1}+1)}$. Then $p(x)$ is a palindromic polynomial of degree $2m = n-3$, and it has positive coefficients by virtue of Lemma \ref{lem:Greene} with $\alpha = \frac{2\pi}{n+1}$, $\theta = \frac{5\pi}{n+1}$. Construct $u \in \mathbb R_+^{n-2}$ and ${\bf u} \subset \mathbb R_+^n$ as above from the coefficients of $p(x)$. Then the angles $\zeta_j = \frac{(2j+3)\pi}{n+1}$, $j = 1,\dots,m$, satisfy conditions (a),(b) of Lemma \ref{lem:circulant_odd}. By Theorem \ref{thm:circ_char_odd} $F_{\bf u}$ is one-dimensional and generated by a circulant extremal copositive matrix with minimal circulant zero support set whose elements are given by \eqref{formula_odd}. By the third assertion of Lemma \ref{explicit_elements} the elements of $A$ are explicitly given by
\begin{equation} \label{example_reg}
A_{ij} = \left\{ \begin{array}{rcl} 2(1+2\cos\frac{\pi}{n+1}\cos\frac{3\pi}{n+1}),&\quad&i = j, \\
-2(\cos\frac{\pi}{n+1}+\cos\frac{3\pi}{n+1}),&&|i-j| \in \{1,n-1\}, \\
1,&&|i-j| \in \{2,n-2\}, \\
0,&&|i-j| \in \{3,\dots,n-3\}, \end{array} \right.
\end{equation}
$i,j = 1,\dots,n$.

\medskip

However, Lemma \ref{lem:Greene} allows also for other choices of regularly spaced angles $\zeta_1,\dots,\zeta_m$ or regularly spaced angles $\zeta_1,\dots,\zeta_{m-1},2\pi - \zeta_m$. The following result guarantees the positivity of the coefficients of $p(x)$ also in the case when only $\zeta_2,\dots,\zeta_m$ (or $\zeta_2,\dots,\zeta_{m-1},2\pi-\zeta_m$) are regularly spaced, and the spacing between $\zeta_1$ and $\zeta_2$ is inferior to the spacing between the other angles.

{\lemma \cite[Corollary 1.1]{BDPW91} Let $p(x)$ be a real polynomial with nonnegative coefficients, and let $x_0$ be the root of $p(x)$ which has the smallest argument among all roots of $p(x)$ in the upper half-plane. If $x_1$ is any number such that $|x_1| \geq |x_0|$ and $\it{Re}\,x_1 \leq \it{Re}\,x_0$, then the coefficients of the polynomial $p(x)\frac{(x-x_1)(x-\bar x_1)}{(x-x_0)(x-\bar x_0)}$ are not smaller than the corresponding coefficients of $p(x)$. }

\medskip

As to the general case, we establish the following conjecture.

{\conjecture Let $n \geq 5$ be an integer, and set $m = \lceil \frac{n}{2} \rceil - 2$. Let $\zeta_1,\dots,\zeta_m \in (0,\pi]$ be an increasing sequence of angles such that the fractional part of $\frac{n\zeta_j}{4\pi}$ is in $(0,\frac12)$ for odd $j$ and in $(\frac12,1)$ for even $j$. Define the polynomial $q(x) = \prod_{j=1}^m(x^2 - 2x\cos\zeta_j +1)$, and set $p(x) = q(x)$ for odd $n$ and $p(x) = (x+1)q(x)$ for even $n$. Then the coefficients of $p(x)$ are all positive if and only if $\zeta_j \in (\frac{(2j+2)\pi}{n},\frac{(2j+4)\pi}{n})$ for all $j = 1,\dots,m$. }

We have verified this conjecture for $n \leq 8$.

\section{Matrices of order 6} \label{sec:n6}

In this section we compute all exceptional extremal copositive matrices $A$ of size $6 \times 6$ which have zeros $u^j$ with $\Supp u^j = I_j$, $j = 1,\dots,6$. We assume without loss of generality that $\diag A = (1,\dots,1)$, because every other such matrix lies in the ${\cal G}_6$-orbit of a matrix normalized in this way.

Let $A \in {\cal C}^6$ be exceptional and extremal, $\diag A = (1,\dots,1)$, and let ${\bf u} = \{ u^1,\dots,u^6 \} \subset \mathbb R_+^6$ be zeros of $A$ satisfying $\Supp u^j = I_j$. By Theorem \ref{thm:reg_deg} the matrix $A$ has non-minimal circulant zero support set, and $\rk A_{I_j} = 2$ for all $j = 1,\dots,6$. Therefore $B = \Lambda(A)$ has rank 2 by Corollary \ref{cor:rank_equal}, and $B = x \otimes x + y \otimes y$ for some linearly independent solutions $x = (x_1,x_2,\dots),y = (y_1,y_2,\dots) \in {\cal L}_{\bf u}$.
Extremality of $A$ implies that the linear span of $\{x,y\}$ does not contain non-zero periodic solutions. Indeed, let $v \in \spa\{x,y\}$ be such a solution. Then $v \otimes v \in {\cal P}_{\bf u}$, and $B \pm \varepsilon \cdot v \otimes v \in {\cal F}_{\bf u}$ for $\varepsilon$ small enough. Hence $A$ can be represented as a non-trivial convex combination of the elements $\Lambda^{-1}(B \pm \varepsilon v \otimes v) \in F_{\bf u}$, contradicting extremality.

Shift-invariance of $B$ implies by virtue of Lemma \ref{lem:inv_unitary} that the monodromy $\mathfrak M$ of the system ${\bf S}_{\bf u}$ acts on $x,y \in {\cal L}_{\bf u}$ by
\[ \begin{pmatrix} {\mathfrak M} x \\ {\mathfrak M} y \end{pmatrix} = M\begin{pmatrix} x \\ y \end{pmatrix},
\]
where $M$ is an orthogonal $2 \times 2$ matrix. If $M$ is a reflection, then $\mathfrak M$ has an eigenvector with eigenvalue 1 in the span of $\{x,y\}$, leading to a contradiction by Lemma \ref{periodic1}. Likewise, $M$ cannot be the identity matrix. Hence $M$ is a rotation by an angle $\zeta \in (0,2\pi)$, i.e.,
\begin{equation} \label{monodromy6}
{\mathfrak M}x = \cos\zeta \cdot x - \sin\zeta \cdot y, \qquad {\mathfrak M}y = \sin\zeta \cdot x + \cos\zeta \cdot y.
\end{equation}

By \eqref{band_relation} and by the normalization adopted above we have $x_k^2 + y_k^2 = 1$ for all $k \geq 1$. Without loss of generality we may assume that $x_1 = 1$, $y_1 = 0$, otherwise we rotate the basis $\{x,y\}$ of $\spa\{x,y\}$ appropriately by redefining $x,y$. Then we have
\[ x_k = \cos\sum_{j=1}^{k-1}(\pi-\varphi_j), \qquad y_k = \sin\sum_{j=1}^{k-1}(\pi-\varphi_j)
\]
for some angles $\varphi_1,\varphi_2,\dots \in [0,2\pi)$, for all $k \geq 1$. By virtue of \eqref{monodromy6} we get
\begin{equation} \label{refgl1}
\sum_{j=k}^{k+5} (\pi-\varphi_j) \equiv -\sum_{j=k}^{k+5}\varphi_j \equiv \zeta \quad \mod\ 2\pi\qquad \forall\ k \geq 1.
\end{equation}
It follows that $\varphi_k \equiv \varphi_{k+6}$ modulo $2\pi$ for all $k \geq 1$, and the sequence $\{\varphi_k\}$ is 6-periodic.

We have
\begin{equation} \label{solution6}
B({\bf e}_t,{\bf e}_s) = x_tx_s + y_ty_s = (-1)^{s-t}\cos\sum_{j=\min(t,s)}^{\max(t,s)-1}\varphi_j,\qquad \forall\ t,s \geq 1.
\end{equation}
Conditions \eqref{linear_relation} reduce to $B({\bf e}_t,{\bf e}_{t+3}) = B({\bf e}_{t+3},{\bf e}_{t+6})$ for all $t \geq 1$, which yields
\begin{equation} \label{eq6}
\cos(\varphi_t+\varphi_{t+1}+\varphi_{t+2}) = \cos(\varphi_{t+3}+\varphi_{t+4}+\varphi_{t+5}) = \cos(\zeta+\varphi_t+\varphi_{t+1}+\varphi_{t+2}) \qquad \forall\ t \geq 1.
\end{equation}
By virtue of $\zeta \not= 0$ modulo $2\pi$ it follows that
\begin{equation} \label{refgl2}
\varphi_t+\varphi_{t+1}+\varphi_{t+2} \equiv -(\zeta+\varphi_t+\varphi_{t+1}+\varphi_{t+2})\quad \mod \ 2\pi,
\end{equation}
or equivalently,
\[ \varphi_t+\varphi_{t+1}+\varphi_{t+2} \equiv -\frac{\zeta}{2}\quad \mod\ \pi \qquad \forall\ t \geq 1.
\]
This yields $\varphi_t \equiv \varphi_{t+3}$ modulo $\pi$, or $\varphi_{t+3} = \varphi_t + \delta_t$ with $\delta_t \in \{-\pi,0,\pi\}$, for all $t \geq 1$. 

By Lemma \ref{lem:second_cond} inequalities \eqref{ineq_relations} hold strictly, which yields
\begin{equation} \label{ineq6}
\cos(\varphi_t+\varphi_{t+1}) > \cos(\varphi_{t+2}+\varphi_{t+3}+\varphi_{t+4}+\varphi_{t+5}) = \cos(\zeta+\varphi_t+\varphi_{t+1}) \qquad \forall\ t \geq 1.
\end{equation}
Equivalently, $\varphi_t+\varphi_{t+1} \in 2\pi l_t + (-\frac{\zeta}{2},\pi-\frac{\zeta}{2})$ for some $l_t \in \{0,1,2\}$, for all $t \geq 1$. Replacing $t$ by $t+3$, we get $\varphi_{t+3}+\varphi_{t+4} = \varphi_t+\varphi_{t+1}+\delta_t+\delta_{t+1} \in 2\pi l_{t+3} + (-\frac{\zeta}{2},\pi-\frac{\zeta}{2})$, and therefore $\delta_t+\delta_{t+1} \equiv 0$ modulo $2\pi$, for all $t \geq 1$. 

Combining \eqref{refgl1} with \eqref{refgl2} we get $\varphi_{t+3}+\varphi_{t+4}+\varphi_{t+5}\equiv\varphi_t+\varphi_{t+1}+\varphi_{t+2}$ modulo $2\pi$ for all $t \geq 1$, which leads to $\delta_t+\delta_{t+1}+\delta_{t+2} \equiv 0$ modulo $2\pi$ for all $t \geq 1$.


Therefore $\delta_t = 0$ and $\varphi_{t+3} = \varphi_t$ for all $t \geq 1$. Set $\sigma = \varphi_1 + \varphi_2 + \varphi_3$. Then \eqref{ineq6} reduces to $\cos(\sigma-\varphi_j) > \cos(\sigma+\varphi_j)$, or equivalently $\sin\sigma\sin\varphi_j > 0$ for $j = 1,2,3$. Therefore $\sin\varphi_j$ and $\sin\sigma$ have the same sign for all $j = 1,2,3$. By possibly replacing the solution $y$ by $-y$, we may assume without loss of generality that $\sin\varphi_1 > 0$ and hence $\varphi_j \in (0,\pi)$, $\pi - \varphi_j \in (0,\pi)$ for all $j = 1,2,3$. If $\sigma \in (2\pi,3\pi)$, then $\sum_{j=1}^3(\pi - \varphi_j) \in (0,\pi)$, and $y_k > 0$ for $k = 2,3,4$. Since also $y_1 = 0$, this contradicts the condition that $y$ is a solution of a linear 3-rd order system with positive coefficients. Therefore $\varphi_1+\varphi_2+\varphi_3 < \pi$. By \eqref{band_relation},\eqref{solution6} the matrix $A$ is then given by
\begin{equation} \label{classif6}
\begin{pmatrix} \scriptstyle 1 & \scriptstyle -\cos\varphi_1 & \scriptstyle  \cos(\varphi_1+\varphi_2) & \scriptstyle  -\cos(\varphi_1+\varphi_2+\varphi_3) & \scriptstyle  \cos(\varphi_2+\varphi_3) & \scriptstyle  -\cos\varphi_3 \\ \scriptstyle -\cos\varphi_1 & \scriptstyle  1 & \scriptstyle  -\cos\varphi_2 & \scriptstyle  \cos(\varphi_2+\varphi_3) & \scriptstyle  -\cos(\varphi_1+\varphi_2+\varphi_3) & \scriptstyle  \cos(\varphi_1+\varphi_3) \\ \scriptstyle \cos(\varphi_1+\varphi_2) & \scriptstyle  -\cos\varphi_2 & \scriptstyle  1 & \scriptstyle  -\cos\varphi_3 & \scriptstyle  \cos(\varphi_1+\varphi_3) & \scriptstyle  -\cos(\varphi_1+\varphi_2+\varphi_3) \\ \scriptstyle -\cos(\varphi_1+\varphi_2+\varphi_3) & \scriptstyle  \cos(\varphi_2+\varphi_3) & \scriptstyle  -\cos\varphi_3 & \scriptstyle  1 & \scriptstyle  -\cos\varphi_1 & \scriptstyle  \cos(\varphi_1+\varphi_2) \\ \scriptstyle \cos(\varphi_2+\varphi_3) & \scriptstyle  -\cos(\varphi_1+\varphi_2+\varphi_3) & \scriptstyle  \cos(\varphi_1+\varphi_3) & \scriptstyle  -\cos\varphi_1 & \scriptstyle  1 & \scriptstyle  -\cos\varphi_2 \\ \scriptstyle -\cos\varphi_3 & \scriptstyle  \cos(\varphi_1+\varphi_3) & \scriptstyle  -\cos(\varphi_1+\varphi_2+\varphi_3) & \scriptstyle  \cos(\varphi_1+\varphi_2) & \scriptstyle  -\cos\varphi_2 & \scriptstyle  1 \end{pmatrix}.
\end{equation}

On the other hand, let $\varphi_1,\varphi_2,\varphi_3 > 0$ such that $\varphi_1+\varphi_2+\varphi_3 < \pi$, and consider the matrix $A$ given by \eqref{classif6}. Let $v^1,\dots,v^6 \in \mathbb R_+^6$ be the columns of the matrix
\begin{equation} \label{V6}
V = \begin{pmatrix} \sin\varphi_2 & 0 & 0 & 0 & \sin\varphi_2 & \sin(\varphi_1+\varphi_3) \\ \sin(\varphi_1+\varphi_2) & \sin\varphi_3 & 0 & 0 & 0 & \sin\varphi_3 \\ \sin\varphi_1 & \sin(\varphi_2+\varphi_3) & \sin\varphi_1 & 0 & 0 & 0 \\ 0 & \sin\varphi_2 & \sin(\varphi_1+\varphi_3) & \sin\varphi_2 & 0 & 0 \\ 0 & 0 & \sin\varphi_3 & \sin(\varphi_1+\varphi_2) & \sin\varphi_3 & 0 \\ 0 & 0 & 0 & \sin\varphi_1 & \sin(\varphi_2+\varphi_3) & \sin\varphi_1 \end{pmatrix},
\end{equation}
and define $u^j = v^j + v^{j+1}$, $j = 1,\dots,5$, $u^6 = v^6 + v^1$. By construction the submatrices $A_{I_j}$ are positive semi-definite and of rank 2, and $(u^j)^TAu^j = 0$, $\Supp u^j = I_j$ for all $j = 1,\dots,6$. Moreover, $(u^j)^TAu^{j+1} > 0$ for all $j = 1,\dots,5$ and $(u^6)^TAu^1 > 0$. Hence $A$ is an exceptional copositive matrix by Theorem \ref{first_cond}, and it is extremal with non-minimal circulant zero support set by Theorem \ref{thm:reg_deg}. By (ii.d) of Theorem \ref{thm:reg_deg} the columns of $V$ are minimal zeros of $A$, and every minimal zero of $A$ is a positive multiple of some column of $V$. We have proven the following result.

{\theorem \label{thm:class6} Let $A \in {\cal C}^6$ be exceptional and extremal with zeros $u^1,\dots,u^n$ satisfying $\Supp u^j = I_j$ for all $j = 1,\dots,6$. Then $A$ is in the ${\cal G}_6$-orbit of some matrix of the form \eqref{classif6}, with $\varphi_1,\varphi_2,\varphi_3 > 0$ and $\varphi_1+\varphi_2+\varphi_3 < \pi$. The minimal zero support set of $A$ is given by $\{\{1,2,3\},\{2,3,4\},\{3,4,5\},\{4,5,6\},\{1,5,6\},\{1,2,6\}\}$.
On the other hand, every matrix $A$ of the form \eqref{classif6} with $\varphi_1,\varphi_2,\varphi_3 > 0$ and $\varphi_1+\varphi_2+\varphi_3 < \pi$ is exceptional and extremal, and every minimal zero of $A$ is proportional to one of the columns of the matrix \eqref{V6}. \qed }

{\remark We do not claim that every extremal exceptional copositive matrix in ${\cal C}^6$ with $\diag A = (1,\dots,1)$ and with this minimal zero support set has to be of the form \eqref{classif6}. }

\section{Conclusions}

In this contribution we considered copositive matrices with zeros $u^1,\dots,u^n \in \mathbb R_+^n$ having supports $\Supp u^j = I_j$. Exceptional copositive matrices with this property exist for all matrix sizes $n \geq 5$ and are of two types, in dependence on whether the zeros $u^j$ are minimal or not. The matrices of each type make up an algebraic submanifold of ${\cal S}^n$, of codimensions $n$ and $2n$, respectively. The prototypes of these matrices are the Hildebrand matrices and the Horn matrix, respectively. Explicit examples of such matrices have been given in \eqref{example_reg} and \eqref{example_deg}, respectively. We show that if the zeros $u^j$ are not minimal, then the corresponding exceptional copositive matrices are extremal. If the zeros $u^j$ are minimal, then the corresponding matrices can be extremal only for odd $n$.

Some open questions within this framework remain:
\begin{itemize}
\item Do non-extremal matrices with minimal circulant zero support set exist for odd matrix size?
\item Do matrices with non-minimal circulant zero support set and having minimal zero support set different from $\{I_1',\dots,I_n'\}$ exist?
\item Which types of extremal copositive matrices can appear on the boundary of the submanifolds of matrices with minimal and non-minimal circulant zero support set, respectively?
\item What is the global topology of these submanifolds?
\end{itemize}

\section*{Acknowledgements}

This research was carried out in the framework of \textsc{Matheon} supported by Einstein Foundation Berlin. The paper has also benefitted from the remarks of anonymous referees.

\bibliographystyle{plain}

\appendix

\section{Auxiliary results}

In this section we collect a few auxiliary results of general nature.

{\lemma \label{completion} Let $A \in {\cal S}^n$ and define the ordered index subsets $I = (1,\dots,n-1)$, $I' = (2,\dots,n)$. Suppose that the principal submatrices $A_I,A_{I'}$ are positive semi-definite, and that there exist vectors $u,v \in \mathbb R^n$ such that $u_1,v_n > 0$, $v_1 = u_n = 0$, and $u^TAu = v^TAv = 0$. Then there exists a unique real number $\delta$ such that $A - \delta \cdot E_{1n}$ is positive semi-definite, and this number has the same sign as the product $u^TAv$. Moreover, the vector $u+v$ is in the kernel of $P(\delta)$. In particular, if all elements of $u,v$ except $v_1,u_n$ are positive, then $P(\delta)$ has an element-wise positive kernel vector. }

\begin{proof}
Consider the matrix-valued function $P(\delta) = A - \delta E_{1n}$. All elements of $P(\delta)$ except the upper right (and correspondingly lower left) corner element coincide with the corresponding elements of $A$. Therefore the posed problem on the unknown $\delta$ can be considered as a positive semi-definite matrix completion problem. Namely, we wish to modify the corner elements of $A$ to make it positive semi-definite. By a standard result on positive semi-definite matrix completions \cite{Grone&Johnson&Sa&Wolkowicz:84} there exists a solution $\delta$ such that $P(\delta) \succeq 0$.

On the other hand, for every solution $\delta$ we have $u^TP(\delta)u = v^TP(\delta)v = 0$ and therefore we must have $P(\delta)u = P(\delta)v = 0$ and $u^TP(\delta)v = 0$. Uniqueness of the solution follows and we have the explicit expression $\delta = \frac{u^TAv}{u^TE_{1n}v}$. The first assertion of the lemma now follows from the strict inequality $u^TE_{1n}v > 0$. The second assertion holds because $u+v$ must be in the kernel of $P(\delta)$ along with $u$ and $v$.
\end{proof}

{\lemma \label{corank1_manifold} Define the set ${\cal M} = \{ M \in \mathbb R^{n \times n} \,|\, \det M = 0 \}$, and let $M \in {\cal M}$ be a matrix of corank 1. Let the vectors $u,v \in \mathbb R^n$ be generators of the left and the right kernel of $M$, respectively. Then the orthogonal complement under the Frobenius scalar product $\langle A,B \rangle = \tr(AB^T)$ of the tangent space to ${\cal M}$ at $M$ is generated by the rank 1 matrix $uv^T$. }

\begin{proof}
Let $M = F^0_L(F^0_R)^T$ be a factorization of $M$, where $F^0_L,F^0_R$ are full column rank $n \times (n-1)$ matrices. The set ${\cal M}$ can be written as $\{ F_LF_R^T \,|\, F_L,F_R \in \mathbb R^{n \times (n-1)} \}$. Therefore the tangent space to ${\cal M}$ at $M$ is given by all matrices of the form $F^0_L\Delta_R^T + \Delta_L (F^0_R)^T$, where $\Delta_L,\Delta_R \in \mathbb R^{n \times (n-1)}$ are arbitrary. Therefore a matrix $H \in \mathbb R^{n \times n}$ is orthogonal to ${\cal M}$ at $M$ if and only if $\langle F^0_L\Delta_R^T + \Delta_L (F^0_R)^T,H \rangle = \tr(\Delta_L (F^0_R)^TH^T + \Delta_R (F^0_L)^TH) = 0$ for all $\Delta_L,\Delta_R \in \mathbb R^{n \times (n-1)}$. Equivalently, $HF^0_R = H^TF^0_L = 0$, or $HM^T = H^TM = 0$. The assertion of the lemma now readily follows.
\end{proof}

{\lemma \label{lem:gradients} Let $u^1,\dots,u^n \in \mathbb R^n$ be non-zero vectors, such that no two of these are proportional and the dimension of the linear span $\spa\{u^1,\dots,u^n\}$ has corank at most 1. Then the rank 1 matrices $u^j(u^j)^T$, $j = 1,\dots,n$, are linearly independent. }

\begin{proof}
Without loss of generality we may assume that $u^1,\dots,u^{n-1}$ are linearly independent. In an appropriate coordinate system these vectors then equal the corresponding canonical basis vectors. Then $u^j(u^j)^T = E_{jj}$ for $j = 1,\dots,n-1$. Hence the rank 1 matrices $u^j(u^j)^T$ are linearly dependent only if $u^n(u^n)^T$ is diagonal with at least two non-zero entries, which leads to a contradiction.
\end{proof}

{\lemma \label{corank2_manifold} Define the set ${\cal M} = \{ S \in {\cal S}^n \,|\, \det S_{(1,\dots,n-1)} = \det S_{(2,\dots,n)} = \det (S_{ij})_{i = 1,\dots,n-1;j = 2,\dots,n} = 0 \}$, and let $S \in {\cal M}$ be a positive semi-definite matrix of corank 2. Suppose there exists a basis $\{u,v\} \subset \mathbb R^n$ of $\ker S$ such that $u_1 \not= 0$, $v_n \not= 0$, $u_n = v_1 = 0$. Then there exists a neighbourhood ${\cal U} \subset {\cal S}^n$ of $S$ such that a matrix $S' \in {\cal U}$ is positive semi-definite of corank 2 if and only if $S' \in {\cal U} \cap {\cal M}$. }

\begin{proof}
By $u_n = 0$ the subvector $u_{(1,\dots,n-1)}$ is in the kernel of $S_{(1,\dots,n-1)}$. Hence $S_{(1,\dots,n-1)}$ is of corank at least 1. If $w \in \ker S_{(1,\dots,n-1)}$ is another kernel vector, then $w' = (w^T, 0)^T \in \mathbb R^n$ is in the kernel of $S$ and must therefore be proportional to $u$, because $v$ cannot be in $\spa\{ u,w \}$ by $v_n \not= 0$. Therefore $S_{(1,\dots,n-1)}$ is positive semi-definite of corank 1. Similarly, the submatrix $S_{(2,\dots,n)}$ is positive semi-definite of corank 1.

Let $S' \in {\cal M}$ be close to $S$. Then by continuity the $n-2$ largest eigenvalues of $S'_{(1,\dots,n-1)}$ are positive, and the remaining eigenvalue is zero by definition of ${\cal M}$. Therefore $S'_{(1,\dots,n-1)} \succeq 0$ and $\rk S' \geq n-2$. The kernel of $S'_{(1,\dots,n-1)}$ is close to that of $S_{(1,\dots,n-1)}$, hence there exists a vector $u'$, close to $u$, such that $(u')^TS'u' = 0$ and $u'_1 \not= 0$, $u'_n = 0$. Similarly, $S'_{(2,\dots,n)}$ is positive semi-definite and there exists a vector $v'$, close to $v$, such that $(v')^TS'v' = 0$ and $v'_n \not= 0$, $v'_1 = 0$.

By $u'_1 \not= 0$ the first column of the submatrix $S'_{(1,\dots,n-1)}$ is a linear combination of the other columns. These $n-2$ columns must therefore be linearly independent. It follows that the submatrix $S'_{(2,\dots,n-1)}$ is positive definite. Therefore the $(n-2) \times (n-1)$ submatrix $(S'_{ij})_{i = 2,\dots,n-1;j = 2,\dots,n}$ has full row rank, and every vector in its right kernel must be proportional to $v'_{(2,\dots,n)}$. Hence the right kernel of the singular submatrix $(S'_{ij})_{i = 1,\dots,n-1;j = 2,\dots,n}$ must also be generated by $v'_{(2,\dots,n)}$. Similarly, the left kernel of this submatrix is generated by $u'_{(1,\dots,n-1)}$, and we get that $(u')^TS'v' = 0$.

By possibly replacing $u'$ by $-u'$ or $v'$ by $-v'$, we may enforce $u'_1 > 0$, $v'_n > 0$. By Lemma \ref{completion} we then have that $S'$ is positive semi-definite. Now both $u'$ and $v'$ are in the kernel of $S'$, and these vectors are linearly independent. Therefore $\rk S' \leq n-2$, and $S'$ is of corank 2.

On the other hand, every matrix $S'$ of corank 2 is in ${\cal M}$. This completes the proof.
\end{proof}

{\lemma \label{lem:Toeplitz_rk1sum} Let $T$ be a singular real symmetric positive semi-definite Toeplitz matrix of rank $k$ and with an element-wise nonnegative non-zero kernel vector $u$. Then there exist distinct angles $\zeta_0 = \pi$, $\zeta_1,\dots,\zeta_m \in (0,\pi)$ and positive numbers $\lambda_0,\dots,\lambda_m$, where $m = \lfloor \frac{k}{2} \rfloor$, such that $T = \sum_{j=0}^m \lambda_j T(\zeta_j)$ for odd $k$ and $T = \sum_{j=1}^m \lambda_j T(\zeta_j)$ for even $k$, where $T(\zeta)$ is the symmetric Toeplitz matrix with first row $(1,\cos\zeta,\cos 2\zeta,\dots)$. Moreover, if $p(x)$ is the polynomial whose coefficients equal the entries of $u$, then $e^{\pm i\zeta_j}$ are roots of $p(x)$ for $j = 1,\dots,m$, and $-1$ is a root of $p(x)$ if $k$ is odd. }

\begin{proof}
Any positive semi-definite Toeplitz matrix $T$ of rank $k$ can be represented as a weighted sum of $k$ rank 1 positive semi-definite Toeplitz matrices with positive weights. Each of these rank 1 matrices is complex Hermitian with first row $(1,e^{i\zeta},e^{2i\zeta},\dots)$ for some $\zeta \in [0,2\pi)$, and the angles $\zeta$ are pairwise distinct. If $T$ is singular, then the weights and the angles are determined uniquely \cite[Chapter 3]{Zarrop}. For real $T$ the complex rank 1 matrices appear in complex conjugate pairs, each of which sums to a Toeplitz matrix of the form $T(\zeta)$ with $\zeta \in (0,\pi)$.

The kernel of $T$ equals the intersection of the kernels of the rank 1 summands. The angle $\zeta = 0$ then cannot appear in the sum due to the presence of an element-wise nonnegative non-zero kernel vector. Hence the angle $\zeta = \pi$, which corresponds to the only remaining real rank 1 Toeplitz matrix, appears in the sum if and only if $k$ is odd.

Finally, $u$ is in the kernel of every rank 1 summand in the decomposition of $T$. This directly yields the last assertion of the lemma.
\end{proof}



{\lemma \label{lambda_eq_aux} Let $n \geq 5,r,s$ be positive integers, $\gamma \in \mathbb R$ arbitrary, and $\zeta_1,\dots,\zeta_s \in [0,\pi]$. Then the system of $r$ linear equations
\begin{eqnarray*}
\sum_{j=1}^s \lambda_j (\cos(n-k)\zeta_j - \cos k\zeta_j) &=& 0,\qquad k = \left\lceil \frac{n}{2} \right\rceil + 1 - r,\dots,\left\lceil \frac{n}{2} \right\rceil - 1, \\
\sum_{j=1}^s \lambda_j \left[\cos\left(n-\left(\left\lceil \frac{n}{2} \right\rceil-r\right)\right)\zeta_j - \cos \left(\left\lceil \frac{n}{2} \right\rceil-r\right)\zeta_j\right] &=& \gamma \nonumber
\end{eqnarray*}
in the unknowns $\lambda_1,\dots,\lambda_s$ is equivalent to the system
\[ \begin{pmatrix} 1 & \cdots & 1 \\ \cos\zeta_1 & \cdots & \cos\zeta_s \\ \vdots & & \vdots \\ \cos^{r-1}\zeta_1 & \cdots & \cos^{r-1}\zeta_s \end{pmatrix} \begin{pmatrix} \lambda_1\sin\frac{\zeta_1}{2}\sin\frac{n\zeta_1}{2} \\ \vdots \\ \lambda_s\sin\frac{\zeta_s}{2}\sin\frac{n\zeta_s}{2} \end{pmatrix} = \begin{pmatrix} 0 \\ \vdots \\ 0 \\ -2^{-r}\gamma \end{pmatrix}
\]
for odd $n$ and
\[ \begin{pmatrix} 1 & \cdots & 1 \\ \cos\zeta_1 & \cdots & \cos\zeta_s \\ \vdots & & \vdots \\ \cos^{r-1}\zeta_1 & \cdots & \cos^{r-1}\zeta_s \end{pmatrix} \begin{pmatrix} \lambda_1\sin\zeta_1\sin\frac{n\zeta_1}{2} \\ \vdots \\ \lambda_s\sin\zeta_s\sin\frac{n\zeta_s}{2} \end{pmatrix} = \begin{pmatrix} 0 \\ \vdots \\ 0 \\ -2^{-r}\gamma \end{pmatrix}
\]
for even $n$. }

\begin{proof}
Due to the identity $\cos(2\alpha-\beta) - \cos\beta = -2\sin\alpha\sin(\alpha-\beta)$ the system is equivalent to
\[ \begin{pmatrix} \sin(\frac{n}{2} - \lceil \frac{n}{2} \rceil + 1)\zeta_1 & \cdots & \sin(\frac{n}{2} - \lceil \frac{n}{2} \rceil + 1)\zeta_s \\ \vdots & & \vdots \\ \sin(\frac{n}{2} - \lceil \frac{n}{2} \rceil + r)\zeta_1 & \cdots & \sin(\frac{n}{2} - \lceil \frac{n}{2} \rceil + r)\zeta_s \end{pmatrix} \begin{pmatrix} \lambda_1\sin\frac{n\zeta_1}{2} \\ \vdots \\ \lambda_s\sin\frac{n\zeta_s}{2} \end{pmatrix} = \begin{pmatrix} 0 \\ \vdots \\ 0 \\ -\frac{\gamma}{2} \end{pmatrix}.
\]
Column $j$ of the coefficient matrix above is given by $(\sin\frac{\zeta_j}{2},\sin\frac{3\zeta_j}{2},\dots,\sin\frac{(2r-1)\zeta_j}{2})^T$ for odd $n$ and by $(\sin\zeta_j,\sin2\zeta_j,\dots,\sin r\zeta_j)^T$ for even $n$. Apply the formula $\sin k\varphi = (e^{(k-1)i\varphi} + e^{(k-3)i\varphi} + \dots + e^{-(k-1)i\varphi})\sin\varphi$ to the coefficients in column $j$ with $\varphi = \frac{\zeta_j}{2}$ for odd $n$ and $\varphi = \zeta_j$ for even $n$.

By adding to each row of the coefficient matrix appropriate multiples of the rows above it, we may obtain a matrix whose column $j$ is given by $\sin\frac{\zeta_j}{2}(1,e^{i\zeta_j}+e^{-i\zeta_j},\dots,(e^{i\zeta_j}+e^{-i\zeta_j})^{r-1})^T$ for odd $n$ and by $\sin\zeta_j(1,e^{i\zeta_j}+e^{-i\zeta_j},\dots,(e^{i\zeta_j}+e^{-i\zeta_j})^{r-1})^T$ for even $n$. The right-hand side of the system does not change under this operation. Recalling that $e^{i\zeta}+e^{-i\zeta} = 2\cos\zeta$, we obtain the equivalent systems in the formulation of the lemma.
\end{proof}

{\corollary \label{lambda_eq} Let $n \geq 5$ be an integer, and set $m = \lceil \frac{n}{2} \rceil - 2$. Let $\zeta_1,\dots,\zeta_m \in (0,\pi]$ be distinct angles. Then the inhomogeneous linear system of the $m$ equations
\begin{eqnarray} \label{linear_inhom}
\sum_{j=1}^m \lambda_j (\cos(n-k)\zeta_j - \cos k\zeta_j) &=& 0,\qquad 3 \leq k \leq m + 1, \\
\sum_{j=1}^m \lambda_j (\cos(n-2)\zeta_j - \cos 2\zeta_j) &=& -1 \nonumber
\end{eqnarray}
in the unknowns $\lambda_1,\dots,\lambda_m$ has a solution if and only if among the angles $\zeta_j$ there are no multiples of $\frac{2\pi}{n}$, and this solution is unique and given by
\begin{equation} \label{solution_lambda}
\lambda_j = \left\{ \begin{array}{rcl}
\left(2^m\sin\frac{\zeta_j}{2}\sin\frac{n\zeta_j}{2}\prod_{l\not=j}(\cos\zeta_j-\cos\zeta_l)\right)^{-1},&\qquad& n\ \mbox{odd}, \\
\left(2^m\sin\zeta_j\sin\frac{n\zeta_j}{2}\prod_{l\not=j}(\cos\zeta_j-\cos\zeta_l)\right)^{-1},&\qquad& n\ \mbox{even}.
\end{array} \right.
\end{equation} }

\begin{proof}
Apply Lemma \ref{lambda_eq_aux} with $r = s = m$. Since the cosine function is strictly monotonous on $(0,\pi]$, the Vandermonde matrix in the equations in this lemma is non-singular. Using explicit formulas for the inverse of the Vandermonde matrix \cite{MaconSpitzbart58}, we obtain for every $j = 1,\dots,m$ that $\lambda_j\sin\frac{\zeta_j}{2}\sin\frac{n\zeta_j}{2} = \frac{1}{2^m\prod_{l \not= j}(\cos\zeta_j - \cos\zeta_l)}$ for odd $n$ and $\lambda_j\sin\zeta_j\sin\frac{n\zeta_j}{2} = \frac{1}{2^m\prod_{l \not= j}(\cos\zeta_j - \cos\zeta_l)}$ for even $n$. The claim of the corollary now easily follows.
\end{proof}

{\lemma \label{coef_positive} Let $n \geq 1$ be an integer and $\zeta_1,\dots,\zeta_m \in I$ a strictly increasing sequence of angles, where $I = (0,\pi]$ if $n$ is odd and $I = (0,\pi)$ if $n$ is even. Then the following are equivalent:
\begin{itemize}
\setlength{\itemsep}{1pt}
\setlength{\parskip}{0pt}
\setlength{\parsep}{0pt}
\item[(i)] the fractional part of $\frac{n\zeta_j}{4\pi}$ is in $(0,\frac12)$ for odd $j$ and in $(\frac12,1)$ for even $j$;
\item[(ii)] the coefficients $\lambda_1,\dots,\lambda_m$ given by \eqref{solution_lambda} are all positive.
\end{itemize} }

\begin{proof}
For every $j = 1,\dots,m$ we have $\sin\zeta_j > 0$ if $n$ is even and $\sin\frac{\zeta_j}{2} > 0$ if $n$ is odd. Further, the cosinus function is strictly decreasing on $[0,\pi]$, and therefore $\sign\prod_{l\not=j}(\cos\zeta_j-\cos\zeta_l) = (-1)^{j-1}$ for $j = 1,\dots,m$. Hence $\lambda_j > 0$ if and only if $\sign\sin\frac{n\zeta_j}{2} = (-1)^{j-1}$, which in turn is equivalent to condition (i).
\end{proof}


\section{Extremality of exceptional circulant matrices} \label{sec:appB}

In this section we provide Lemma \ref{lem:circ_extremal}, which is the technically most difficult part of the proof of Lemma \ref{lem:circulant_odd}. Its proof uses a bit more advanced mathematical concepts, in particular, linear group representations.

Let $L_C \subset {\cal S}^n$ be the subspace of circulant matrices, i.e., matrices which are invariant with respect to simultaneous circular shifts of the row and column indices. Let $P_S: {\cal S}^n \to {\cal S}^n$ be the linear map which corresponds to a circular shift by one entry. Then $P_S^n$ is the identity map $\it{Id}$, and $P_S$ generates a symmetry group which is isomorphic to the cyclic group $C_n$. However, symmetric circulant matrices possess yet another symmetry. Namely, they are {\it persymmetric}, i.e., invariant with respect to reflection about the main skew diagonal. Denote this reflection by $P_P: {\cal S}^n \to {\cal S}^n$. The symmetry group generated by $P_S$ and $P_P$ is isomorphic to the dihedral group $D_n$, which has $2n$ elements. The action of $D_n$ on ${\cal S}^n$ defines a linear unitary representation $R_{{\cal S}^n}$ of $D_n$ of dimension $\frac{n(n+1)}{2}$, which decomposes into a direct sum of irreducible representations. This decomposition corresponds to an orthogonal decomposition of ${\cal S}^n$ into a direct sum of invariant subspaces.

In this section we suppose $n \geq 5$ to be an odd number. Then the group $D_n$ has two irreducible representations of dimension 1. Both of these send $P_S$ to the identity $\it{Id}$. The element $P_P$ is sent to $\it{Id}$ by the trivial representation and to $-\it{Id}$ by the other 1-dimensional representation. The $\frac{n+1}{2}$-dimensional subspace $L_C$ is exactly the invariant subspace corresponding to the trivial representation. The other 1-dimensional representation does not participate in $R_{{\cal S}^n}$, because there is no non-zero matrix $A \in {\cal S}^n$ such that $P_S(A) = A$ and $P_P(A) = -A$. Besides the 1-dimensional representations, $D_n$ possesses $\frac{n-1}{2}$ 2-dimensional representations $R_1,\dots,R_{\frac{n-1}{2}}$. Here $R_k$ sends $P_S$ to the rotation of the 2-dimensional representing subspace by the angle $\frac{2\pi k}{n}$, and $P_P$ to a reflection of this subspace.

Set $m = \frac{n-3}{2}$ and let $\zeta_1,\dots,\zeta_m \in (0,\pi)$ be distinct angles in increasing order. Define the polynomial $p(x) = \prod_{j=1}^m (x^2 - 2x\cos\zeta_j + 1)$, which has roots $e^{\pm i\zeta_j}$, and let $u \in \mathbb R^{n-2}$ be its coefficient vector. Note that $p$ is palindromic, i.e., $u$ does not change if the order of its entries is inverted. For $j = 1,\dots,n$, define vectors $u^j \in \mathbb R^n$ such that $u^j_{I_j} = u$ and $u^j_i = 0$ for all $i \not\in I_j$. We do not make further assumptions on the collection ${\bf u} = \{ u^1,\dots,u^n \}$, in particular, we do not demand $u$ to be positive.

For every $k = 1,\dots,\frac{n-1}{2}$ we also define the polynomial
\[ p_k(x) = e^{-2im\pi k(n-1)/n}p(e^{i\pi k(n-1)/n}x) = \Pi_{j=1}^m (x - e^{i(\zeta_j - \pi k(n-1)/n)})(x-e^{i(-\zeta_j-\pi k(n-1)/n)})
\]
and denote by $v^k \in \mathbb R^{n-2}$ its coefficient vector, with $v^k_1 = 1$ being the coefficient of $x^{2m}$. Then the roots of $p_k$ all lie on the unit circle and equal $e^{i(\pm\zeta_j-\pi k(n-1)/n)}$, $j = 1,\dots,m$. The elements of $v^k$ are given by $v^k_l = e^{-\pi ik(l-1)(n-1)/n}u_l$, which by the oddity of $n$ equals $e^{\pi ik(n+1-l)(n-1)/n}u_l$, $l = 1,\dots,n-2$.

For given $\bf u$, defined by angles $\zeta_1,\dots,\zeta_m$ as above, we are interested in the linear subspace $L_{\bf u} \subset {\cal S}^n$ of symmetric matrices $A$ such that $A_{I_j}u = 0$ for all $j = 1,\dots,n$. If $A \in L_{\bf u}$ is a solution of the corresponding linear system of equations, then the matrices $P_S^k(A)$ obtained from $A$ by circular shifts of the row and column indices are also solutions. Moreover, since $u$ is palindromic, $P_P(A)$ is a solution too. Therefore $L_{\bf u}$ is an invariant subspace under the action of the group $D_n$, and this action defines a linear representation $R_{\bf u}$ of $D_n$ on $L_{\bf u}$. This representation decomposes into a direct sum of irreducible representations and induces an orthogonal decomposition of $L_{\bf u}$ into a direct sum of invariant subspaces $L_{{\bf u},\it{Id}},L_{{\bf u},k}$, $k = 1,\dots,\frac{n-1}{2}$, corresponding to the trivial irreducible representation and the representations $R_k$ of $D_n$. As noticed above, the non-trivial 1-dimensional representation of $D_n$ does not contribute to $R_{{\cal S}^n}$ and hence neither to $R_{\bf u}$.

The subspace $L_{{\bf u},\it{Id}}$ then consists exactly of those matrices $A$ which are circulant and satisfy $A_{I_1}u = 0$. The next result characterizes the subspaces $L_{{\bf u},k}$.

{\lemma \label{lem:Lk} Let $A \in L_{{\bf u},k}$ be a non-zero matrix. Then there exists a complex symmetric matrix $B$ such that $\it{Re} B, \it{Im} B \in L_{{\bf u},k}$, $A \in \spa\{ \it{Re} B, \it{Im} B \}$, and $B$ can be represented as a Hadamard product $B = C \circ H$ of a non-zero real symmetric circulant matrix $C$ and a Hankel rank 1 matrix $H$ given element-wise by $H_{jl} = e^{\pi i k(n+1-j-l)(n-1)/n}$. Moreover, we have $C_{I_1}v^k = 0$. }

\begin{proof}
Let $L$ be a 2-dimensional subspace which is invariant under the action of $D_n$ and such that $A \in L \subset L_{{\bf u},k}$. Then the action of $D_n$ on $L$ is given by the representation $R_k$, i.e., there exists a basis $\{ B_1,B_2 \}$ of $L$ such that
\[ \begin{pmatrix} P_S(B_1) \\ P_S(B_2) \end{pmatrix} = \begin{pmatrix} \cos\frac{2\pi k}{n} & \sin\frac{2\pi k}{n} \\ -\sin\frac{2\pi k}{n} & \cos\frac{2\pi k}{n} \end{pmatrix} \begin{pmatrix} B_1 \\ B_2 \end{pmatrix}, \qquad \begin{pmatrix} P_P(B_1) \\ P_P(B_2) \end{pmatrix} = \begin{pmatrix} 1 & 0 \\ 0 & -1 \end{pmatrix} \begin{pmatrix} B_1 \\ B_2 \end{pmatrix}.
\]
Setting $B = B_1 + iB_2$, we obtain that $P_S(B) = e^{-2\pi ik/n}B$, $P_P(B) = \overline{B}$.

The first condition is equivalent to the condition $B_{jl} = e^{2\pi ik/n}B_{j_-l_-}$ for all $j,l = 1,\dots,n$, where $j_-,l_- \in \{1,\dots,n\}$ are the unique indices satisfying $j_- \equiv j-1$ and $l_- \equiv l-1$ modulo $n$. For every $j,l = 1,\dots,n$ we then have
\[ B_{jl} = \left\{ \begin{array}{rcl} \exp\left(\frac{2\pi ik}{n}\frac{j+l-n-1}{2}\right) \cdot B_{j-\frac{j+l-n-1}{2},l-\frac{j+l-n-1}{2}},&\qquad& (j+l)\ \mbox{even}; \\
\exp\left(\frac{2\pi ik}{n}\frac{j+l-1}{2}\right) \cdot B_{j+n-\frac{j+l-1}{2},l-\frac{j+l-1}{2}},&\qquad& j < l,\ (j+l)\ \mbox{odd}; \\
\exp\left(\frac{2\pi ik}{n}\frac{j+l-1}{2}\right) \cdot B_{j-\frac{j+l-1}{2},l+n-\frac{j+l-1}{2}},&\qquad& j > l,\ (j+l)\ \mbox{odd}.
\end{array} \right.
\]
Since $n$ is odd, the first factors in this representation of the elements of $B$ form the Hankel matrix $H$, while the second factors form a circulant matrix $C$ which coincides with $B$ on the main skew-diagonal. However, the condition $P_P(B) = \overline{B}$ implies that the elements on the main skew diagonal of $B$ are real. Hence $C$ is also real. This yields the claimed representation $B = C \circ H$. Since $B \not= 0$, we also have $C \not= 0$.

Since $B_1,B_2 \in L_{\bf u}$, we have that $B_{I_1}u = 0$. However, for all $l = 1,\dots,n-2$ we have $(B_{I_1}u)_l = e^{-\pi ikl(n-1)/n}(C_{I_1}v^k)_l$ by definition of $v^k$, and hence also $C_{I_1}v^k = 0$. This completes the proof.
\end{proof}

Recall that for a circulant matrix $A$, the submatrix $A_{I_1}$ is Toeplitz. The next result shows that the conditions $Tu = 0$ or $Tv^k = 0$ impose stringent constraints on a real symmetric Toeplitz matrix $T$.

{\lemma \label{lem:Toeplitz_linear_comb} Let $\varphi_1,\dots,\varphi_d \in (-\pi,\pi]$ be distinct angles, and let $w \in \mathbb R^{d+1}$ be the coefficient vector of the polynomial $p(x) = \prod_{k=1}^d (x - e^{i\varphi_k})$, with $w_1 = 1$ being the coefficient of $x^d$. Let $\Xi \subset [0,\pi]$ be the set of angles $\xi$ such that either both $\pm\xi$ appear among the angles $\varphi_k$, $k = 1,\dots,d$, or $\xi = \varphi_j = \pi$ for some $j$. Let $T$ be a real symmetric Toeplitz matrix satisfying $Tw = 0$. Then $T \in \spa \{T_{\xi} \}_{\xi \in \Xi}$, where $T_{\xi}$ is the symmetric Toeplitz matrix with first row $(1,\cos\xi,\dots,\cos d\xi)$, and the linear span is over $\mathbb R$. }

\begin{proof}
Let $\varphi_{d+1},\dots,\varphi_{2d+1} \in (-\pi,\pi]$ be such that $\varphi_j$ are mutually distinct for $j = 1,\dots,2d+1$. For any $\varphi \in (-\pi,\pi]$, let $h_{\varphi} = (e^{ik\varphi})_{k = 0,\dots,d} \in \mathbb C^{d+1}$ be a column vector and $H_{\varphi} = h_{\varphi}h^*_{\varphi}$ the corresponding complex Hermitian rank 1 Toeplitz matrix. Then $H_{\varphi_1},\dots,H_{\varphi_{2d+1}}$ are linearly independent over the complex numbers, and hence also over the reals. Indeed, the matrix $H_{\varphi}$ contains $2d+1$ distinct elements $e^{-id\varphi},\dots,e^{id\varphi}$. However, the vectors $(e^{-id\varphi_j},\dots,e^{id\varphi_j})$, $j = 1,\dots,2d+1$, are linearly independent, because suitable multiples of these vectors can be arranged into a Vandermonde matrix. Therefore $H_{\varphi_1},\dots,H_{\varphi_{2d+1}}$ form a basis of the $(2d+1)$-dimensional real vector space of complex Hermitian $(d+1) \times (d+1)$ Toeplitz matrices.

The real symmetric Toeplitz matrix $T$ is also complex Hermitian and can hence be written as a linear combination $T = \sum_{j=1}^{2d+1} \alpha_j H_{\varphi_j}$, $\alpha_j \in \mathbb R$. Now $h^*_{\varphi_j}w = 0$ for $j = 1,\dots,d$ by construction of $w$, and therefore $Tw = \sum_{j=d+1}^{2d+1} \alpha_j (h^*_{\varphi_j}w) h_{\varphi_j} = 0$. But the vectors $h_{\varphi_{d+1}},\dots,h_{\varphi_{2d+1}}$ again form a Vandermonde matrix and are hence linearly independent. Moreover, we have $h^*_{\varphi_j}w = e^{-id\varphi_j}p(e^{i\varphi_j}) \not= 0$ for all $j = d+1,\dots,2d+1$. It follows that $\alpha_j = 0$ for all $j = d+1,\dots,2d+1$, and $T = \sum_{j=1}^d \alpha_j H_{\varphi_j}$.

Now the first column of the imaginary part $\it{Im} H_{\varphi_j}$ is given by $\it{Im} h_{\varphi_j} = (0,\sin\varphi_j,\dots,\sin d\varphi_j)^T$, and hence $\sum_{j=1}^d \alpha_j \sin l\varphi_j = 0$ for all $l = 1,\dots,d$. As in the proof of Lemma \ref{lambda_eq_aux}, we may use the formula $\sin l\varphi = \sin\varphi(e^{(l-1)i\varphi} + e^{(l-3)i\varphi} + \dots + e^{-(l-1)i\varphi})$ to rewrite this system of linear equations in the coefficients $\alpha_j$ as
\[ \begin{pmatrix} 1 & \cdots & 1 \\ \cos\varphi_1 & \cdots & \cos\varphi_d \\ \vdots & & \vdots \\ \cos^{d-1}\varphi_1 & \cdots & \cos^{d-1}\varphi_d \end{pmatrix} \begin{pmatrix} \alpha_1\sin\varphi_1 \\ \vdots \\ \alpha_d\sin\varphi_d \end{pmatrix} = 0.
\]
It follows that a coefficient $\alpha_j$ can only be non-zero if either $\sin\varphi_j = 0$, in which case $\varphi_j \in \Xi$ and $H_{\varphi_j} = T_{\varphi_j}$, or there exists another index $j' \in \{1,\dots,d\}$ such that $\varphi_{j'} = -\varphi_j$ and $\alpha_{j'} = \alpha_j$, and therefore $|\varphi_j| \in \Xi$ and $\alpha_j H_{\varphi_j} + \alpha_{j'}H_{\varphi_{j'}} = 2\alpha_jT_{|\varphi_j|}$. This completes the proof.
\end{proof}

The restrictions on the Toeplitz matrices translate into the following restrictions on matrices $A \in L_{{\bf u},\it{Id}}$ and on the circulant factor in Lemma \ref{lem:Lk}.

{\corollary \label{cor:circ_solution} Let $n \geq 5$ be an odd integer, set $m = \frac{n-3}{2}$, and let $\varphi_1,\dots,\varphi_{2m} \in (-\pi,\pi]$ be distinct angles, such that there are no multiples of $\frac{2\pi}{n}$ among them. Define the vector $w \in \mathbb R^{n-2}$ and the set $\Xi \subset [0,\pi]$ as in Lemma \ref{lem:Toeplitz_linear_comb}, with $d = 2m$.

Let $L_w \subset L_C$ be the linear subspace of real symmetric circulant matrices $C$ such that $C_{I_1}w = 0$. Then $\dim L_w \leq 1$. If $\dim L_w = 1$, then $\Xi$ has $m$ elements. In particular, for $\dim L_w = 1$ either the values $e^{i\varphi_j}$ group into $m$ complex-conjugate pairs, or they group into $m-1$ complex conjugate pairs and among the two remaining values one equals $-1$. }

\begin{proof}
Since $\varphi_j \not= 0$ for all $j = 1,\dots,2m$, the set $\Xi$ can have at most $m$ elements. Set $r = |\Xi|$.

Let $C \in L_w$ be non-zero. By Lemma \ref{lem:Toeplitz_linear_comb} the Toeplitz matrix $T = C_{I_1}$ can be written as a non-zero linear combination of $\{ T_{\xi} \}_{\xi \in \Xi}$, $T = \sum_{j=1}^r \lambda_jT_{\xi_j}$, $\xi_j \in \Xi$ for all $j = 1,\dots,r$. Since the elements of $C_{I_1}$ determine the circulant matrix $C$ completely by Lemma \ref{lem:irred_N}, we also have $\dim L_w \leq r$. Since $C \not= 0$, the set $\Xi$ contains at least one element.

For $n = 5$ we then get $r = m = 1$, which proves the assertion in this case.

Suppose that $n \geq 7$. We have $C_{1k} = C_{1,n+2-k}$ for all $k = 2,\dots,n$. It follows that $T_{1k} = T_{1,n+2-k}$ for all $k = 4,\dots,\frac{n+1}{2}$. This yields the linear homogeneous system of equations $\sum_{j=1}^r \lambda_j\cos k\xi_j = \sum_{j=1}^r \lambda_j\cos(n-k)\xi_j$, $k = 3,\dots,\frac{n-1}{2}$, in the coefficients $\lambda_j$. By Lemma \ref{lambda_eq_aux} this system is equivalent to the system
\[ \begin{pmatrix} 1 & \cdots & 1 \\ \cos\xi_1 & \cdots & \cos\xi_r \\ \vdots & & \vdots \\ \cos^{m-2}\xi_1 & \cdots & \cos^{m-2}\xi_r \end{pmatrix} \begin{pmatrix} \lambda_1\sin\frac{\xi_1}{2}\sin\frac{n\xi_1}{2} \\ \vdots \\ \lambda_r\sin\frac{\xi_r}{2}\sin\frac{n\xi_r}{2} \end{pmatrix} = 0.
\]
Since there are no multiples of $\frac{2\pi}{n}$ among the angles $\varphi_j$, there are no multiples of $\frac{2\pi}{n}$ among the $\xi_j$ neither, and hence $\sin\frac{\xi_j}{2}\sin\frac{n\xi_j}{2} \not= 0$ for all $j = 1,\dots,r$. It follows that the coefficient matrix in the system above has a non-trivial kernel, implying $r > m-1$ and hence $r = m$.

However, the coefficient matrix has full row rank $m-1$, and therefore the solution $\lambda_1,\dots,\lambda_m$ is determined by the angles $\xi_1,\dots,\xi_m$ up to multiplication by a common scalar. Hence $C_{I_1}$ and by Lemma \ref{lem:irred_N} also $C$ is determined up to a scalar factor, and $\dim L_w \leq 1$.
\end{proof}

We are now in a position to prove the result we need for Lemma \ref{lem:circulant_odd}.

{\lemma \label{lem:circ_extremal} Let $n \geq 5$ be odd. Set $m = \frac{n-3}{2}$ and let $\zeta_1,\dots,\zeta_m \in (0,\pi)$ be distinct angles in increasing order such that the fractional part of $\frac{n\zeta_j}{4\pi}$ is in $(0,\frac12)$ for odd $j$ and in $(\frac12,1)$ for even $j$. Let $u \in \mathbb R^{n-2}$ be the coefficient vector of the polynomial $p(x) = \prod_{j=1}^m (x^2 - 2x\cos\zeta_j + 1)$ and let $L_{\bf u} \subset {\cal S}^n$ be the subspace of all matrices $A$ satisfying $A_{I_j}u = 0$ for all $j = 1,\dots,n$. Then the condition $\dim L_{\bf u} > 1$ implies that $\zeta_j = \frac{(2j-1)\pi}{n}$ for all $j = 1,\dots,m$. In particular, the vector $u$ has negative elements. }

\begin{proof}
By assumption there are no multiples of $\frac{2\pi}{n}$ among the angles $\pm\zeta_j$ and, since $n-1$ is even, also among the angles $\pm\zeta_j-\frac{\pi k(n-1)}{n}$ for all $j = 1,\dots,m$ and $k = 1,\dots,\frac{n-1}{2}$. Note also that $\frac{\pi k(n-1)}{n}$ is not a multiple of $2\pi$ for these values of $k$, and hence the angles $\pm\zeta_j-\frac{\pi k(n-1)}{n}$ are obtained from the angles $\pm\zeta_j$ by a shift by a non-zero value modulo $2\pi$.

The solution space $L_{\bf u}$ is a direct sum of the subspaces $L_{{\bf u},\it{Id}}$ and $L_{{\bf u},k}$, $k = 1,\dots,\frac{n-1}{2}$ corresponding to the irreducible representations of the group $D_n$. By Corollary \ref{cor:circ_solution} with $w = u$ we have $\dim L_{{\bf u},\it{Id}} \leq 1$. Therefore $\dim L_{\bf u} > 1$ implies that $L_{{\bf u},k}$ is non-zero for some $k$.

Let $A \in L_{{\bf u},k}$ be a non-zero matrix. Then by Lemma \ref{lem:Lk} there exists a non-zero real symmetric circulant matrix $C$ satisfying $C_{I_1}v^k = 0$. By Corollary \ref{cor:circ_solution} with $w = v^k$ the roots $e^{i(\pm\zeta_j-\pi k(n-1)/n)}$, $j = 1,\dots,m$, of $p_k(x)$ either group into $m$ complex-conjugate pairs, or they group into $m-1$ complex conjugate pairs and among the two remaining roots one equals $-1$. We shall now show that this condition uniquely determines the angles $\zeta_j$.

For $l = 1,\dots,n$, define open arcs $a_l = \{ e^{i\varphi} \,|\, \varphi \in (\frac{2\pi(l-1)}{n},\frac{2\pi l}{n}) \}$ of length $\frac{2\pi}{n}$ on the unit circle. Then $e^{\pm i\zeta_j},e^{i(\pm\zeta_j-\pi k(n-1)/n)} \in \bigcup_{l=1}^n a_l$ for all $j = 1,\dots,m$, $k = 1,\dots,\frac{n-1}{2}$. Moreover, by the assumptions on $\zeta_j$ we can have $e^{i\zeta_j} \in a_l$ only if the parity condition $l \equiv j$ modulo 2 holds. Since $\zeta_1,\dots,\zeta_m$ is an increasing sequence, we also have that each interval $a_l$ contains at most one of the values $e^{i\zeta_1},\dots,e^{i\zeta_m}$. We consider two cases.

1. $\zeta_m > \frac{(n-1)\pi}{n}$. Then $e^{\pm i\zeta_m} \in a_{(n+1)/2}$. Any other arc $a_l$ contains at most one of the values $e^{\pm i\zeta_j}$. Hence there are exactly 4 of these arcs containing no such value, and these are located symmetrically about the real axis. Now consider how the rotated values $e^{i(\pm\zeta_j-\pi k(n-1)/n)}$, $j = 1,\dots,m$, are distributed over the arcs $a_l$. In a similar way there must be 4 arcs, call them $\alpha_1,\alpha_2,\alpha_3,\alpha_4$, containing no value and one arc, call it $\beta$, containing two values, but $\beta \not= a_{(n+1)/2}$. We distinguish again two possibilities.

1.1. If the complex conjugate arc to $\beta$ is one of the arcs $\alpha_1,\dots,\alpha_4$, then the two complex values contained in $\beta$ are not matched by complex conjugate values.

1.2. If none of the arcs $\alpha_1,\dots,\alpha_4$ is complex conjugate to $\beta$, then at least one of these arcs, let it be $\alpha_1$, is matched by a complex conjugate arc containing exactly one value. Hence at least one of the values in $\beta$ and the value in the complex conjugate arc to $\alpha_1$ are not matched by complex conjugate values.

It follows in both cases that there are at least two complex values among $e^{i(\pm\zeta_j-\pi k(n-1)/n)}$ which are not matched by a complex conjugate, leading to a contradiction.

2. $\zeta_m < \frac{(n-1)\pi}{n}$. Since $e^{i\zeta_m} \not\in a_{m+1}$ by the parity condition, we must have $e^{i\zeta_j} \in a_j$, $e^{-i\zeta_j} \in a_{n+1-j}$ for all $j = 1,\dots,m$. Now we again consider the distribution of the values $e^{i(\pm\zeta_j-\pi k(n-1)/n)}$, $j = 1,\dots,m$. There are 3 consecutively located arcs which do not contain any of the values $e^{i(\pm\zeta_j-\pi k(n-1)/n)}$, call them $\alpha_1,\alpha_2,\alpha_3$. The remaining arcs contain exactly one value each. The arcs $\alpha_1,\alpha_2,\alpha_3$ are obtained from the arcs $a_{(n-1)/2},a_{(n+1)/2},a_{(n+3)/2}$ by multiplication with $e^{-i\pi k(n-1)/n}$ and are not located symmetrically about the real axis. Hence at least one value among $e^{i(\pm\zeta_j-\pi k(n-1)/n)}$, $j = 1,\dots,m$, is not matched by a complex conjugate. It follows that among these values the value $-1$ must appear, and two of the arcs $\alpha_1,\alpha_2,\alpha_3$ must be complex conjugate to each other. The second condition is only possible if these two arcs border the point $z = 1$ on the unit circle. Equivalently, $e^{i\pi k(n-1)/n}$ must lie on the boundary of the arc $a_{(n+1)/2}$, implying $k = 1$. Recall that the value among $e^{i(\pm\zeta_j-\pi(n-1)/n)}$, $j = 1,\dots,m$ which lies in the interval $a_{(n+1)/2}$ must equal $-1$. This is the value $e^{i(-\zeta_1-\pi(n-1)/n)}$, and hence $\zeta_1 = \frac{\pi}{n}$. Finally, using that the values $e^{i(\zeta_j-\pi(n-1)/n)}$ and $e^{i(-\zeta_{j+1}-\pi(n-1)/n)}$ are mutually complex conjugate for all $j = 1,\dots,m-1$, we obtain the values $\zeta_j = \frac{(2j-1)\pi}{n}$ for all $j = 1,\dots,m$.

Thus we have $p(x) = \prod_{j=1}^m (x^2 - 2x\cos\frac{(2j-1)\pi}{n} + 1)$. The coefficient of the linear term is given by $u_1 = -2\sum_{j=1}^m \cos\frac{(2j-1)\pi}{n} < 0$, which concludes the proof.
\end{proof}

\end{document}